\documentclass{amsart}
\usepackage{amssymb,amsxtra,enumerate}

\newcommand{\Cnh}{\mathbb{C}^{\hat n}}
\newcommand{\C}{\mathbb{C}}
\newcommand{\bC}{\mathbb{C}}
\newcommand{\CN}{\mathbb{C}^N}
\newcommand{\CNplus}{\mathbb{C}^{N+1}}
\newcommand{\CNh}{\mathbb{C}^{\hat N}}

\newcommand{\N}{\mathbb{N}}
\newcommand{\da}{\delta}
\newcommand{\Da}{\Delta}
\newcommand{\io}{\iota}
\newcommand{\al}{\alpha}
\newcommand{\eps}{\epsilon}
\newcommand{\bt}{\beta}
\newcommand{\gm}{\gamma}
\newcommand{\Gm}{\Gamma}
\newcommand{\sg}{\sigma}

\newcommand{\om}{\omega}
\newcommand{\ld}{\lambda}
\newcommand{\gth}{\theta}

\newcommand{\td}{\tilde}

\newcommand{\vardop}[2]{\frac{\partial #1}{\partial #2}}
\newcommand{\br}[1]{\langle#1 \rangle}

\newcommand{\p}{\prime}
\newcommand{\crb}[1]{\mathcal{V} (#1)}
\newcommand{\proj}{\pi}
\newcommand{\GLk}{{GL(\C^k)}}
\newcommand{\CRder}[1]{D_{#1}}
\newcommand{\atlas}{\mathcal{A}}
\newcommand{\fatA}{\mathbf{A}}
\newcommand{\dblind}{\blacktriangle}
\newcommand{\indexset}{\mathcal{K}}

\newcommand{\indexsettwo}{\mathcal{K}_{\hat n}}
\newcommand{\stringA}{\mathbb{A}}
\newcommand{\Span}{\text{\rm span}}
\newlength{\extendaxesby}

\DeclareMathOperator{\spanc}{span}

\DeclareMathOperator{\pr}{proj}
\DeclareMathOperator{\id}{id}
\newtheorem{thm}{Theorem}
\newtheorem{lem}[thm]{Lemma}

\newtheorem{cor}[thm]{Corollary}

\theoremstyle{definition}

\newtheorem{rem}{Remark}

\begin{document}
\title[Finite jet determination of CR embeddings]
{Finite jet determination of constantly degenerate CR embeddings}%
\author{Peter Ebenfelt}%
\address{Department of Mathematics, University of California at
San Diego, La Jolla, CA 92093}%
\email{pebenfel@math.ucsd.edu}%
\author{Bernhard Lamel}
\address{Department of Mathematics, University of Illinois,
Urbana}
\email{lamelb@member.ams.org}%
\thanks{The first
  author is a Royal Swedish Academy of Sciences Research Fellow
  supported by a grant from the Knut and Alice Wallenberg
Foundation and also by DMS-0100110. The second author was
supported by the ANACOGA research network.}
\begin{abstract}
We prove finite jet determination results for smooth CR
embeddings which are of constant degeneracy, using the method of
complete systems. As an application, we derive a reflection
principle for mappings between a Levi-nondegenerate hypersurface
in $ \CN$ and  a Levi-nondegenerate hypersurface in $\CNplus$. We
also give an independent proof of the reflection principle
for mappings between strictly pseudoconvex hypersurfaces
in any codimension due to Forstneric \cite{FO1}.
\end{abstract}
\maketitle

\section{Introduction} In this paper, we shall derive a complete
system of differential equations for CR mappings between a smooth
hypersurface in $\CN$ and a smooth hypersurface in $\CNh$ (see
e.g.\ Theorem \ref {T:complsys}). This system of differential
equations will then be used to prove finite determination results
(Theorems \ref{T:codimension1}, \ref{T:strpcx}, and
\ref{T:nondegenerate}) for such mappings. For historical
information and background on the finite determination problem,
we refer the reader to the introduction of \cite{E5} or the
survey article \cite{BERbull}. As an application of our main
results, we also prove a new reflection principle for smooth
mappings between real-analytic Levi nondegenerate hypersurfaces
(Corollary \ref{C:reflection1}).

To formulate our results precisely, we shall need to introduce
some notation. Let $M$ be a real hypersurface in $\bC^N$, $N\geq
2$, and let us denote by $\mathcal V\subset \bC TM$ its CR bundle;
the reader is referred e.g.\ to \cite{BERbook} for basic notions
and facts about CR manifolds. Recall that a mapping $f\colon M\to
\bC^k$ is called CR if $f_*(\mathcal V)\subset T^{0,1}\bC^k$,
where $f_*$ denotes the induced mapping of the tangent spaces (the
push forward). This is equivalent to saying that, near every point
$p\in M$, there is a neighborhood $U\subset M$ of $p$ in which
\begin{equation}\label{eq-CR}
L_{\bar{a}}f_j=0,\quad a=1,\ldots,n,\ j=1,\ldots k,
\end{equation}
where $n=N-1$ and $L_{\bar 1},\ldots,L_{\bar n}$ is a basis for
the CR vector fields in $U$. We use here the convention that e.g.\
$L_{\bar{a}}=\overline{L_a}$ so that the $L_1,\ldots, L_n$ form a
basis for the vector fields valued in $\overline{\mathcal V}$ near
$p$ on $M$.

Suppose that $\hat M\subset \bC^{\hat N}$ is a real hypersurface
and $f\colon M\to \bC^{\hat N}$ is a CR mapping sending $M$ into
$\hat M$. Let $p$ be a  point in $M$ and $\hat\rho$ a local
defining function for $\hat M$ near $\hat p:=f(p)\in \hat M$. The
mapping $f$ is called {\it transversal} at $p$ if $f_*(\bC T_pM)$
is not contained in the subspace (of codimension 1) $\hat
{\mathcal V}_{\hat p}+\overline{\hat {\mathcal V}}_{\hat
p}\subset\bC T_{\hat p}\hat M$. One can show (see section
\ref{SS:pullback}) that if $f$ is transversal and $M$ is Levi
nondegenerate at $p\in M$, then $f$ is in fact a local embedding
near $p$.

Following Lamel \cite{L1}, we define an increasing sequence of
subspaces $E_k(q)\subset\bC^{\hat N}$, for $q\in M$ near $p$,
associated to the germ of the mapping $f\colon (M,q)\to \hat M$ as
follows. Let $L_{\bar 1},\ldots, L_{\bar n}$ be a basis for the CR
vector fields on $M$ near $p$ as above and define
\begin{equation}\label{eq-ksdeg1}
E_k(q):=\Span\{L^{\bar J}(\hat \rho_{Z'}\circ f)(q) : J\in
(\mathbb Z_+)^n, |J|\leq k\} \subset \bC^{\hat N},
\end{equation}
where $\hat \rho_Z=\big({\partial\hat \rho}/{\partial
Z'_j}\big)_{1\leq j\leq \hat N}$ in some local coordinate system
$Z'$ near $\hat p$; we use here standard multi-index notation $
L^{\bar J}:=L_{\bar 1}^{\bar J_1}\ldots L_{\bar n}^{\bar J_n}$ and
$|J|=J_1+\ldots+J_n$. One can show (cf.\ \cite{L1}) that the
increasing sequence of numbers $d_k(q):=\dim E_k(q)$ is
independent of the choice of local defining function $\hat \rho$
and coordinates $Z'$, as well as of the choice of basis of the CR
vector fields $L_{\bar 1},\ldots, L_{\bar n}$. We shall say, again
following Lamel ({\it loc.\ cit.}) that $f\colon (M,q)\to \hat M$
is $(k_0,s)$-degenerate at $q$ if $s$ is the minimum of the
decreasing sequence of codimensions of $E_k$, i.e.\
$$s:=\min_k(\hat N-d_k),$$ and $k_0$ is the smallest integer such that this
minimum is attained. The case where $s=0$ is actually a
nondegenerate situation, and we say that $f$ is
$k_0$-nondegenerate if it is $(k_0,0)$-degenerate. If $f$ is
$(k_0,s)$-degenerate at a point $q$ and the degeneracy $s$ is
constant in a nieghborhood of $q$, then we shall say that $f$ is
constantly $(k_0,s)$-degenerate at $q$. The notion of
$(k_0,s)$-degeneracy will be presented in a more intrinsic form in
section \ref{SS: conditions}.

The following are the main results in this paper.

\begin{thm}\label{T:codimension1}
Let $M\subset\CN$ and $\hat M \subset \CNplus$ be smooth
hypersurfaces which are Levi-nondegenerate at $p_0\in M$ and $\hat
p_0\in \hat M$, and $f:M\to\hat M$ a smooth CR mapping with
$f(p_0)= \hat p_0$ which is constantly $(k_0,s)$-degenerate and
transversal at $p_0$. Then $f$ is uniquely determined by its
$2k_0$-jet $j^{2k_0}_{p_0} f$ at $p_0$ in a neighbourhood of
$p_0$. That is, there exists a neighbourhood $U$ of $p_0$ in $M$
such that if $g: M\to \hat M$ is another smooth CR map of constant
degeneracy $s$ with $j^{2k_0}_{p_0} f = j^{2k_0}_{p_0} g$, then
$f|_{U} = g|_U$.
\end{thm}

\begin{thm}\label{T:strpcx}
Let $M\subset\CN$ and $\hat M \subset \CNh$ be smooth
hypersurfaces, $p_0\in M$, $\hat p_0\in \hat M$, $\hat M$ strictly
pseudoconvex at $\hat p_0$, and $f:M\to\hat M$ a smooth CR
mapping with $f(p_0)= \hat p_0$ which is constantly
$(k_0,s)$-degenerate at $p_0$. Then $f$ is uniquely determined by
its $2k_0$-jet $j^{2k_0}_{p_0} f$ at $p_0$ in a neighbourhood of
$p_0$. That is, there exists a neighbourhood $U$ of $p_0$ in $M$
such that if $g: M\to \hat M$ is another smooth CR map of constant
degeneracy $s$ with $j^{2k_0}_{p_0} f = j^{2k_0}_{p_0} g$, then
$f|_{U} = g|_U$.
\end{thm}

Let us note here that in the setting of Theorem~\ref{T:strpcx},
the mapping $f$ is automatically transversal (see e.g. the paper
\cite{FO1} by Forstneric). The reader should also notice that in
Theorem \ref{T:strpcx}, where the hypersurfaces are required to be
strictly pseudoconvex, there is no restriction on the codimension
$\hat N-N$ of the mapping $f$ whereas in Theorem
\ref{T:codimension1}, where no pseudoconvexity is assumed, the
codimension is restricted to be one. The reason for this is easily
seen by considering the case where $M\subset \bC^N$ is given by
the equation
$$
\Im z_{1}=\sum_{k=2}^{N}|z_k|^2
$$
and $\hat M\subset \bC^{N+2}$ by
$$
\Im z_{1}=\sum_{k=2}^{N}|z_k|^2+|z_{N+1}|^2-|z_{N+2}|^2.
$$
The mapping
\begin{equation}\label{E:h}
(z_1,\ldots,z_N)\mapsto (z_1,\ldots,z_{N},f(z),-f(z)),
\end{equation}
where $f$ is an arbitrary holomorphic function with $f(0)=0$,
sends $(M,0)$ into $(\hat M,0)$. It is not difficult to see that
if $k_0$ is the smallest integer $\geq 2$ for which there is a
multi-index $\alpha\in\mathbb Z_+^{N-1}$ with
$\partial^\alpha_{z'} f(0)\neq 0$, where $z'=(z_2,\ldots,z_N)$,
then the mapping \eqref{E:h} is constantly $(k_0,1)$-degenerate at
$0$. Clearly, you can find two different such functions $f_1,f_2$
whose jets at $0$ agree up to arbitrary order and, hence, Theorem
\ref{T:codimension1} is false in codimension $2$.

We shall see that in fact, given a mapping $f$ as in the preceding
theorems, there exists a dense, open subset $M_f$ of $M$ such
that for each $p\in M_f$, $f$ is constantly $(k_0,s)$- degenerate
for some $s$ and furthermore, that $k_0 \leq \hat N -N+1 - s\leq
\hat N - N + 1$ (see Lemma~\ref{L:bounds}). Hence, we have the
following corollaries.

\begin{cor}\label{C:codimension1}
Let $M\subset\CN$ and $\hat M \subset \CNplus$ be smooth
hypersurfaces which are Levi-nondegenerate,
and $f:M\to\hat M$ and $g: M\to \hat M$  smooth, transversal CR
mappings.
 If for any $p_0$ in the dense open subset
$M_f\cap M_g$ of $M$ we have that $j^{2(\hat N - N+1)}_{p_0} f =
j^{2(\hat N - N+1)}_{p_0} g$, then
$f = g$.
\end{cor}

\begin{cor}\label{C:strpcx}
Let $M\subset\CN$ and $\hat M \subset \CNh$ be smooth strictly
pseudoconvex hypersurfaces and $f:M\to\hat M$ and $g: M \to \hat
M$ smooth CR mappings. If for any $p_0$ in the dense open subset
$M_f\cap M_g$ of $M$ we have that $j^{2(\hat
N - N+1)}_{p_0} f = j^{2(\hat N - N+1)}_{p_0} g$, then $f = g$.
\end{cor}

As an application of Theorem \ref{T:codimension1}, we obtain the
following reflection principle.

\begin{cor}\label{C:reflection1}
Let $M\subset\CN$ and $\hat M \subset \CNplus$ be real-analytic
hypersurfaces which are Levi-nondegenerate, and $f:M\to\hat M$ \ a
smooth transversal CR mapping. Then there exists an open dense
subset $M_0\subset M$ such that for any $p\in M_0$, there exists
an open neighbourhood $V$ of $p$ in $\CN$ and a holomorphic map
$F\colon V\to \CNh$ such that $F|_{V\cap M} = f$.
\end{cor}

Let us give the (short) proof here. We first point out that given
a mapping $f$ as in Corollary \ref{C:reflection1}, there exists a
dense, open subset $M_0$ of $M$ such that for each $p\in M_0$, $f$
is constantly $(k_0,s)$-degenerate, for some $s$, in an open
neighborhood of $p$ (see Lemma \ref{L:global}); in fact, one may
choose $M_0$ such that $f$ is locally constantly
$(k_0,s)$-degenerate on $M_0$ with $s\leq \hat N-N$ and $k_0 \leq
\hat N -N+1 - s\leq \hat N - N + 1$ (see Lemma~\ref{L:bounds}).
For points $p_0\in M_0$, we can use a result of the second author
(\cite{L2}, Theorem 6) to conclude that the formal Taylor series
of $f$ at $p_0$ converges (in some neighbourhood V of $p_0$) to a
holomorphic mapping $F: V \to \CNh$ with $F(M\cap V)\subset \hat
M$. By Theorem~\ref{T:codimension1}, $F|_M = f$, which proves the
corollary.

The same argument (using Lemma~\ref{L:global} and Theorem
\ref{T:strpcx}) also gives an independent proof of the following
theorem due to Forstneric \cite{FO1}.

\begin{cor}\label{C:reflection2}
Let $M\subset\CN$ and $\hat M \subset \CNplus$ be real-analytic
strictly pseudoconvex hypersurfaces and $f:M\to\hat M$ a smooth CR
mapping. Then there exists an open dense subset $M_0\subset M$
such that for any $p\in M_0$, there exists an open neighbourhood
$V$ of $p$ in $\CN$ and a holomorphic map $F\colon V\to \CNh$ such
that $F|_{V\cap M} = f$.
\end{cor}

A special case for which we also derive finite determination
results  without any conditions on $\hat N$ is the case of {\em
nondegenerate} mappings. These have already been investigated by
Kim and Zaitsev \cite{KZ1} in the general codimension case.
However, our method is a bit different, and we have an easy,
independent proof of the following result for nondegenerate
mappings between hypersurfaces.

\begin{thm}\label{T:nondegenerate}
Let $M\subset\CN$ and $\hat M \subset \CNh$ be smooth
hypersurfaces, $p_0\in M$ a point of finite type, and let
$f:M\to\hat M$ be a smooth CR mapping which is $k_0$-nondegenerate
at $p_0$. Then $f$ is uniquely determined by its $2 k_0$-jet
$j^{2k_0}_{p_0} f$ at $p_0$ in a neighbourhood of $p_0$. That is,
there exists a neighbourhood $U$ of $p_0$ in $M$ such that if $g:
M\to \hat M$ is a smooth CR map with $j^{2k_0}_{p_0} f =
j^{2k_0}_{p_0} g$, then $f|_{U} = g|_U$.
\end{thm}

In order to prove Theorems \ref{T:codimension1}, \ref{T:strpcx},
and \ref{T:nondegenerate} we will construct complete systems of
differential equations for all the cases considered. To formulate
the result more precisely, we shall use the following notation.
For an open subset $U\subset M$, we shall denote by $J^{m} ( U,
\hat M )$ the space of $m$-jets of smooth (not necessarily CR)
mappings $F\colon U\to \hat M$ and by $j^m_p F$ the $m$-jet  of
such a mapping at $p\in U$ (which we shall think of as the
collection of all derivatives of $F$, in some local coordinate
system near $p$, up to order $m$). The result needed to prove
Theorems \ref{T:codimension1} and \ref{T:strpcx} can now be
formulated as follows.

\begin{thm}\label{T:complsys} Let $M\subset\CN$, $\hat M \subset \CNh$ be
smooth
hypersurfaces, $p_0\in M$, $\hat p_0\in\hat M$, $f: M\to \hat M$ a
transversal smooth CR map with $f(p_0)= \hat p_0$ of constant
degeneracy $(k_0, s)$ at $p_0$, and assume that either $M\subset
\CN$ and $\hat M\subset \C^{N+1}$ are Levi nondegenerate at $p_0$
and $\hat p_0$ respectively, or that $M$ and $\hat M$ are strictly
pseudoconvex at $p_0$ and $\hat p_0$ respectively (with no
restrictions on $\hat N$). Then there exists a neighbourhood $U$
of $p_0$ in $M$ and a smooth function $\phi : U \times J^{2k_0 +
2} ( U, \hat M )\to  J^{2k_0 + 2} ( U, \hat M )$ which only
depends on $M$ and $\hat M$ such that
\begin{equation}\label{E:complsys}
  j^{2k_0 + 3}_x f = \phi (x, j^{2k_0+2}_x f), \quad x\in U.
\end{equation}
Furthermore, there exists a neigbourhood $V$ of $j^{2k_0 +2}_{p_0}
f$ such that if $g: M\to \hat M$ a transversal smooth CR map with
$g(p_0)= \hat p_0$ of constant degeneracy $(k_0, s)$ at $p_0$ with
$j^{2 k_0 + 2}_{p_0} g \in V$, then (for some neighbourhood $U^\p$
of $p_0$)
\begin{equation}\label{E:complsys2}
  j^{2k_0 + 3}_x g = \phi (x, j^{2k_0+2}_x g), \quad x\in U^\p.
\end{equation}
\end{thm}

Theorems \ref{T:codimension1} and \ref{T:strpcx} follow from
Theorem \ref{T:complsys} using standard arguments about uniqueness
of solutions of complete systems of differential equations (see
e.g.\ \cite{BCG}) and the fact that the $m$-jet of any mapping
$f\colon M\to \hat M$, as in Theorems \ref{T:codimension1} and
\ref{T:strpcx}, at $p_0$, for any $m$ (and hence in particular the
$2k_0+2$-jet) is determined by its $2k_0$ jet at $p_0$ (\cite{L2},
Theorem 8). Thus, to prove Theorems \ref{T:codimension1} and
\ref{T:strpcx} it suffices to prove Theorem \ref{T:complsys}.
Also, by a remark given earlier, we could drop the assumption of
transversality in the strictly pseudoconvex case, since it is
satisfied automatically.

The organization of the paper is as follows: In section
\ref{S:CRvbundles}, we review the definition of a CR vector
bundle, give definitions of the nondegeneracy conditions
introduced above in this context, and introduce some essential
notation and tools. Following that, we give the proof of Theorem
\ref{T:nondegenerate}. In the next two sections, we derive some
jet reflection identities which we use in section
\ref{S:completesystem} to prove Theorem~\ref{T:complsys}.

\section{CR vector bundles}\label{S:CRvbundles} In this section, we
define the notion of a CR vector bundle. This notion is not new
and has been extensively used in the literature, but we give a
thorough definition nonetheless. We also give some examples, and
introduce the canonical derivative operator $ \CRder{}$ (which is
the same as $\bar
\partial_b$ in most of the literature). For our purposes, the
local formalism of this derivative operator is more important than
the global formalism developed alongside it.
\subsection{Definition} Let $M$ be a CR manifold. Let us first recall some
basic definitions. We are given a subbundle $\crb{M} \subset \C T
M$ (which we will refer to as the {\em CR bundle} of $M$) which
satisfies
\begin{equation}\label{E:CRbundle}
  [ \crb{M} , \crb{M} ] \subset \crb{M}, \quad \text{and } \crb{M}
  \cap \overline{\crb{M}} = \{ 0 \}.
\end{equation}
In \eqref{E:CRbundle}, $[ \crb{M} , \crb{M} ] \subset \crb{M}$
means that for all sections $X,Y$ of $\crb{M}$, the section
$[X,Y]$ of $\C T M$ takes values in $\crb{M}$, and $0$ stands for
the $0$ section of $\C T M$. If $M$ is a CR manifold, we also say
that $M$ has a CR structure, and sometimes we even refer to
$\crb{M}$ as a CR structure on $M$. If $U\subset M$ is open, a
section $\bar L \in \Gm(U,\crb{M})$ is called a CR vector field
over $U$. A CR function on $U$ is one which is annihilated by all
CR vector fields on $U$. If $M$, $\hat M$ are two CR manifolds
and $f:M \to \hat M$ is a map, say of class $C^1$, we say that
$f$ is CR provided $f_* (\crb{M}) \subset \crb{\hat M}$. Let us
also recall that if $X$ is a complex manifold, then $X$ has a
natural CR structure given by $\crb{X} = T^{(0,1)} X$.

We consider a complex vector bundle $E$ over $M$ with projection
$\proj:E\to M$. A vector bundle atlas $\atlas$ consists of pairs
$(U_\al,\phi_\al)$ where $\{U_\al\colon \al\in A\}$  is an open
cover of $M$, and $\phi_\al : \proj^{-1} (U_\al)\to  U_\al \times
\C^k$ is a diffeomorphism such that $\proj|_{\proj^{-1}(U_\al)} =
\pr_1 \circ \phi_\al$, where $\pr_l$ denotes the
projection on the $l$th component in a product. The
chart change $\phi_\al \circ \phi_\bt^{-1}\colon (U_\al \cap
U_\bt) \times \C^k\to (U_\al \cap U_\bt) \times \C^k$,
is of the form $\phi_\al \circ \phi_\bt^{-1} (x,v) = (x,\phi_{\al
\bt} (x) v)$, where the {\em transition functions} $\phi_{\al
\bt}$ are smooth functions (on $U_\al\cap U_\bt$) valued in $\GLk$
satisfying the usual cocyle conditions (cf.\ e.g\ \cite{GH}).

We say that $(E,\atlas)$ is a {\em CR vector bundle} if for $\al,
\bt \in A$  the transition functions $\phi_{\al \bt}$  are CR. We
refer to $\atlas$ as a {\em CR atlas} for $E$, and if a CR atlas
is fixed, we simply say that $E$ is a CR vector bundle.

A vector bundle chart $(U, \phi)$  for $E$ is said to be {\em
compatible} with $\atlas$ if $\phi_\al \circ \phi^{-1} (x,v) = (x,
\psi_\al (x) v)$ where $\psi_\al : U\cap U_\al\to \GLk$ is CR for
all $\al \in A$; we say that $(U,\phi)$ is a {\em CR chart} (for
$(E,\atlas)$), and that $U$ is a {\em trivializing neighbourhood}.
If $\atlas$ has the property that every CR chart $(U,\phi)$ for
$(E,\atlas)$ belongs to $\atlas$, we say that $\atlas$ is a {\em
maximal CR atlas} (for E). By including all compatible charts we
can associate to $\atlas$ a unique maximal atlas $\atlas_{max}$
with the property that $\atlas \subset \atlas_{max}$. Two CR
vector bundles $(E,\atlas)$ and $(E,\atlas^\p)$  are compatible if
$\atlas_{max} = \atlas_{max}^\p$. Equivalently, the CR vector
bundles $(E,\atlas)$ and $(E,\atlas^\p)$ are compatible if every
chart $(U_\al,\phi_\al)\in \atlas$ is a CR chart for
$(E,\atlas^\p)$ and vice versa. If a CR atlas $\atlas$ for $E$ is
fixed, we drop the atlas from the notation, understanding that if
we refer to $E$ we really refer to $(E,\atlas)$.

We should point out that a given vector bundle $E$ may have many
incompatible CR vector bundle structures. For instance, the
trivial line bundle $E=M\times \C$ with the trivial CR atlas
$\atlas$ consisting of the single chart $(M,\phi)$ with
$\phi(x,u)=(x,u)$ is compatible with the CR atlas $\atlas'$
consisting of the single chart $(M,\phi')$ with
$\phi'(x,u)=(x,f(x)u)$, where $f\colon M\to \C$ is some
nonvanishing function, if and only if $f$ is CR.

\subsection{The canonical CR structure} In what follows, $(E,\atlas)$
will denote a CR vector bundle over the CR manifold $M$. We shall
keep the CR atlas $\atlas$ fixed and simply refer to the CR vector
bundle $(E,\atlas)$ as $E$. In this section, we shall define a CR
structure on $E$. The main point is that the CR structures induced
on $\proj^{-1}(U_\al)$ by the product CR structures on $U_\al
\times \C^k$, via  the diffeomorphisms $\phi_\al$, fit together
since the transition functions are CR. The CR structure induced on
$E$ in this way can also be characterized by a property of
sections of $E$.

\begin{lem}\label{L:CRstruc} Assume that at each point $p\in M$
there exists a CR function $g$ with ${\bar L} \bar g (p) \neq 0$
for some CR vector field ${\bar L}$. Then there exists a unique CR
structure on $E$ such that a section $\sg\in\Gm(U,E)$ on an open
subset $U\subset M$ is CR if and only if for each CR chart
$(V,\phi)$ of $E$ the components of $\pr_2\circ\phi \circ \sg$ are
CR functions on $U\cap V$. Furthermore, relative to this CR
structure on $E$, the projection $\proj$ is CR, and for each $p\in
M$, there exist a neighbourhood $U$ of $p$ and a local basis of CR
sections for $E$ over $U$.
\end{lem}

Before we start with the proof, let us explain some of the notions
used in the lemma. $M$ is of finite type if the Lie algebra
generated by the CR vector fields and their conjugates, evaluated
at $p\in M$, is all of $\C T_p M$.  A CR section on $U$ is a
section $\sg \in \Gm(U,E)$ which is CR as a map $U\to E$, i.e.
$\sg_* \crb{M }_p \subset \crb{E}_{\sg(p)}$ for all $p\in U$. A
local basis of CR sections (on $U$) is a collection of CR sections
$\sg_1,\dots ,\sg_k\in \Gm(U,E)$ such that $\sg_1 (p), \dots,
\sg_k (p)$ is a basis of $E_p$ for $p\in U$.

We will use that a bundle chart $(U_\al,\phi_\al)$ induces a
splitting of the tangent space $T_P E$ for $P\in \pi^{-1}
(U_\al)$: If we write $\phi_\al(P) = (p,v)$, then $\C T_P E \cong
\C T_{(p,v)} (U_\al \times \C^k) \cong \C T_p M \oplus \C T_v
\C^k$. The former isomorphism is given by the push forward
$(\phi_\al)_*$. To make the latter isomorphism explicit, let
$\io_1$ denote the natural inclusion $\C T_p M\hookrightarrow \C
T_p M \oplus \C T_v \C^k$ given by $\iota_1(X)=X\oplus 0$, and
similarly $\io_2: \C T_v \C^k \hookrightarrow \C T_p M \oplus \C
T_v \C^k $. A tangent vector $\io_1 X + \io_2 Y\in\C T_p M \oplus
T_v \C^k$ with $X\in \C T_{p_0} M$ and $ Y \in \C T_{v_0} \C^k$
acts on the germ of a smooth function $f(p,v)$ at $(p_0,v_0)$ by
$(\io_1 X + \io_2 Y) f = X(f(p,v_0)) + Y(f(p_0,v))$.

\begin{proof}[Proof of Lemma $\ref{L:CRstruc}$]  We define $\crb{E}_P =
(\phi_\al^{-1})_* (\crb{M}_p \oplus T^{(0,1)}_v \C^k)$. We will
show that $\crb{E}$ defined this way is the unique CR structure
claimed in Lemma~\ref{L:CRstruc}. To start, let us first show that
this definition is independent of the chart used. Consider a chart
change $\phi = \phi_\al \circ \phi_\bt^{-1}$. We denote the fiber
coordinate in the chart $\phi_\al$ by $u=(u^1,\dots,u^k)$, the
fiber coordinate in the chart $\phi_\bt$ by $v= (v^1,\dots ,v^k)$,
and let $P_0\in \proj^{-1}(U_\al \cap U_\bt)$, $\phi_\al (P_0) =
(p_0, u_0)$, $\phi_\bt (P_0) = (p_0, v_0)$; thus,
$$
(p,u)=\phi(p,v)=(p,\phi_{\al\bt}(p) v). $$ We must show that
$\phi_*(\crb{M}_{p_0} \oplus T^{(0,1)}_{v_0} \C^k)=\crb{M}_{p_0}
\oplus T^{(0,1)}_{v_0} \C^k$, or equivalently
$\phi_*(\io_1\crb{M}_{p_0} )\subset\crb{M}_p \oplus
T^{(0,1)}_{v_0} \C^k$ and $\phi_*(\io_2T^{(0,1)}_{v_0}
\C^k)\subset\crb{M}_{p_0} \oplus T^{(0,1)}_{v_0} \C^k$. For this
purpose, let $\bar L_{p_0}\in \crb{M}_{p_0}$. Then
\begin{align*}(\phi_{*} {\io_1 \bar L}_{p_0}) f(p,u)
& = {\bar L}_{p_0} f(p,\phi_{\al\bt}(p) v_0) \\
& = {\bar L}_{p_0} f(p, \phi_{\al\bt} (p_0) v_0) + \vardop{f}{u}
{\bar L}_{p_0} (\phi_{\al\bt} (p)) v_0 + \vardop{f}{\bar u} {\bar
L}_{p_0} ( \overline{
\phi_{\al\bt} (p)} ) \bar v_0 \\
&= \io_1 {\bar L}_{p_0} f +  \left( \io_2 \vardop{}{\bar
u}\biggr|_{u_0} f \right)  {\bar L}_{p_0} ( \overline{
\phi_{\al\bt} (p)} ) \bar v_0;
\end{align*}
the last equality holds since $\phi_{\al \bt}$ is CR. Also,
\begin{align*}
(\phi_{*} \io_2 \vardop{}{\bar v^j}\biggr|_{v_0} ) f(p,u) & =
    \io_2 \vardop{}{\bar v^j} \biggr|_{v_0} f(p, \phi_{\al\bt} (p) v) \\
    &= \vardop{}{\bar v^j} \biggr|_{v_0} f(p_0, \phi_{\al\bt} (p_0) v)
    \\
    &= \left( \io_2 \vardop{}{\bar u} \biggr|_{u_0} f \right) \phi_{\al\bt}
(p_0),
\end{align*}
so that $\phi_{*} {\io_1 \bar L}_{p_0} \in \crb{M}_{p_0} \oplus
T^{(0,1)}_{u_0} \C^k$ and $\phi_{*} \io_2 \vardop{}{\bar
v^j}\bigr|_{v_0}\in T^{(0,1)}_{u_0} \C^k$. Hence $\phi_{*}
(\crb{M}_{p_0} \oplus T^{(0,1)}_{v_0} \C^k) = \crb{M}_{p_0} \oplus
T^{(0,1)}_{u_0} \C^k$, as claimed.

It is easy to check that $[\crb{E}, \crb{E}] \subset \crb{E}$
(since the derivatives in the fiber directions and the
derivatives in the base directions commute). It also follows
immediately that $\crb{E} \cap \overline{\crb{E}} = \{ 0 \}$. We
conclude that $\crb{E}$  defines a CR structure on $E$;
furthermore, $\proj_{*} (\crb{E}) = \crb{M}$, so that $\proj$ is
CR as claimed.

We next verify that $\crb{E}$ has the properties claimed in
Lemma~\ref{L:CRstruc}. Let $\sg: U\to E|_U$ be a section, and
$(V,\phi)$ a CR chart with $U\cap V$ nonempty. By the construction
of the CR structure on $E$, $\sigma$ is CR if and only if
$\phi\circ\sigma$ is CR. Since $\sigma$ is a section, we have
$\phi\circ\sigma(p)=(p,\tilde\sigma(p))$, where $\td \sg=
\pr_2\circ\phi \circ \sg = (\sg^1, \dots , \sg^k)^t$. We compute
$(\phi\circ\sigma)_{*} {\bar L}_{p_0}$, where, for ease of
notation, we drop the inclusions $\io_1, \io_2$:
\begin{align*}
((\phi\circ\sigma)_{*} {\bar L}_{p_0}) f(p,u) &= {\bar L}_{p_0} f(p, \td
\sg(p)) \\
 &= {\bar L}_{p_0} f(p, \td \sg(p_0))\, + \, \vardop{f(p_0,u)}{\bar u}
  \Biggr|_{\td \sg(p_0)} {\bar L}_{p_0} \bar { \td \sg }\, +
\,\vardop{f(p_0,u)}{u}
 \Biggr|_{\td \sg(p_0)} {\bar L}_{p_0} \td \sg.
\end{align*}
From this we see that $\phi\circ\sigma$ (and hence $\sg$) is CR if
and only if the components of $\td \sg$ are CR functions.

Now let $p\in M$ and assume that there exists a CR function $g$
with the property that $\bar L \bar g (p) \neq 0$. Let us show
that the CR structure constructed above is unique with the
properties in Lemma \ref{L:CRstruc}. Pick $P_0\in E$ and a CR
chart $(V,\phi)$ with $P_0\in \pi^{-1}(V)$. Write
$\phi(P_0)=(p_0,u_0)$. There is no loss of generality in assuming
that $g(p_0)=u_0$. Using the computation above, first with $\tilde
\sigma=u_0$ and then with $\td \sg = g e_j$ for standard unit
vectors $e_j$, $j=1,\dots,k$, we conclude, since
$\sigma=\phi^{-1}(p,\tilde\sigma)$ is CR for each of the choices
of $\tilde \sigma$ above, that $\crb{E}$ must be the bundle
$(\phi^{-1})_{*}(\crb{M}_{p_0} \oplus T^{(0,1)}_{u_0} \C^k)$. This
completes the proof of the Lemma.
\end{proof}

\begin{rem} {\rm In what follows, unless specifically stated
otherwise, every CR vector bundle will be equipped with the
canonical CR structure provided by Lemma \ref{L:CRstruc}.}
\end{rem}

\subsection{Constructions with CR vector bundles} All standard vector
bundle constructions work in the category of CR vector bundles.
In particular, the sum, tensor product, and exterior product of
CR vector bundles carry the structure of a CR vector bundle, and
the dual bundle $E^*$ is a CR vector bundle. Also, if $\hat E$ is
a CR vector bundle over another CR manifold $\hat M$, and $f:M\to
\hat M$ is a CR map, the pullback bundle $f^{\sharp} \hat E$
inherits the structure of a CR vector bundle.

Let us discuss the construction of the pullback bundle in detail.
We let $E_x = {\hat E}_{f(x)}$ and define a vector bundle
structure on $E= \cup_{x\in M} E_x$ as follows: If $\hat U \subset
\hat M$ is open and $\hat \phi: (\hat \proj)^{-1} (\hat U) \to
\hat U \times \C^k$ is a CR chart for $\hat E$ over $\hat U$,
$(\hat \phi^{-1} \circ (f\times \id))^{-1}$ is a bundle chart for
$E$ over $U=f^{-1}(\hat U)$. Since the transition functions for
$E$ are just the pullbacks (by $f$) of the transition functions
for $\hat E$, $E$ is a CR vector bundle. The map $f:M\to \hat M$
lifts to a differentiable vector bundle map $F:E\to \hat E$. The
canonical CR structure on $\hat E$ is connected with the
structure on $E$ by $\crb{E}_p \subset \crb{\hat E}_{f(p)}$, so
that $F:E\to \hat E$ is a CR map. Also note that the pullback (by
$f$) $f^\sharp \sg\in \Gm(f^{-1}(\hat U), E)$ of a CR section
$\sg\in\Gm(\hat U , \hat E)$ is a CR section of $E$.

\subsection{The derivative operator} In analogy with complex analytic
vector bundles, we can differentiate sections of a CR vector
bundle $E$ (equipped with the canonical CR structure provided by
Lemma \ref{L:CRstruc}) with respect to CR vector fields. To see
this, let $\sg_1,\dots , \sg_k$ be a local basis of CR sections
over some open set $U$. If $\sg \in \Gm(U,E)$ is given by $\sg =
\sum_{j=1}^{k} a^j \sg_j$ with $a^j\in C^\infty (U)$ and ${\bar
L}$ is a CR vector field on $U$, we define
\[
\CRder{{\bar L}} \sg = \sum_{j=1}^{k} (\bar L a^j) \sg_j.
\]
This definition is invariant under changes of the local CR basis
$\sg_j$: Observe that a section $\sg \in \Gm(U,E)$, $\sg =
\sum_{j=1}^{k} a^j \sg_j$ is CR if and only if all the $a^j$ are
CR functions on $U$. Hence, if $\td \sg_1,\dots,\td \sg_k$ is
another basis of CR sections, $\sg_j = \sum_l c_j^{\;\, l}
\td\sg_l$ where the $c_j^{\;\, l}$ are CR functions on $U$. It
follows that $\sg = \sum_l (\sum_j c_j^{\;\, l} a^j)
\td\sg_l=\sum_l\tilde a^l\tilde\sigma_l$, and
\begin{align*}
\sum_l(\bar L\tilde a^l)\tilde\sigma_l  &= \sum_l {\bar L} (\sum_j
c_j^{\;\,l} a^j)
 \td\sg_l \\& =
 \sum_l  (\sum_j c_j^{\;\,l} \bar L a^j) \td\sg_l \\
 &= \sum_j (\bar L a^j) \sg_j,
\end{align*}
so that $\CRder{{\bar L}} \sg$ is independent of the local basis
of CR sections used. This allows us to define $\CRder{{\bar L}}
\sg$ for a section $\sg \in \Gm(U,E)$ on any open set $U\in M$.
Note that $\sg$ is CR if and only if $\CRder{{\bar L}} \sg = 0$
(the $0$ section) for all CR vector fields ${\bar L}$.

We also note that the operator $\om \mapsto \CRder{{\bar L}}\om$
is local, and since $\CRder{a\bar L} \om = a \CRder{{\bar L}}\om$
for any function $a$, $(\CRder{{\bar L}}\om)_p$ depends only on
${\bar L}_p$. Moreover, $\CRder{}$ satisfies the product rule:
\begin{equation}\label{E:productrule}
  \CRder{\bar L} (a\om) = (\bar L a) \om + a \CRder{\bar L} \om.
\end{equation}
If $\omega=(\om_1,\dots ,\om_k)$ is any local (not necessarily CR)
basis of sections over $U$ and ${\bar L}_1,\dots, {\bar L}_n$ is a
local basis of CR vector fields over $U$, we can write
$\CRder{{\bar L}_i} \om_j = \sum_l a_{ij}^{\;\;\, l} \om_l$, (or,
using matrix notation, $\CRder{{\bar L}_i} \om = A_i \om$) and
$A_i=(a_{i j}^{\;\;\,l})_i$ determines $\CRder{{\bar L}} \sigma$
for all sections $\sigma$ over $U$. If we change basis by $\td
\om_j = \sum_m C^{\;m}_j \om_m$, which we will write in matrix
notation as $\td \om = C \om$, we obtain
\begin{align*} \CRder{{\bar L}_i} \td \om
&= \CRder{{\bar L}_i}  C  \om \\
&= ({\bar L}_i C) \om + C (\CRder{{\bar L}_i} \om) \\
&= (({\bar L}_i C) C^{-1} + C A_i C^{-1}) \td \om.
\end{align*}
This means that under a change of basis, the $A_i$ transform as
follows:
\[ \td \om = C \om \Longrightarrow \td A_i = ({\bar L}_i C) C^{-1} + CA_i
C^{-1}.
\]
We can also think of $\CRder{}$ as an operator from sections of
$E$ to sections of $E\otimes \crb{M}^* \simeq E\otimes (\C T^* M /
\mathcal V(M)^\perp)$.  Since $\CRder{}$ satisfies the product
rule \eqref{E:productrule}, we may consider $\CRder{}$ as a
partial connection in $E$.

If $E$ is the trivial bundle $M\times \C$ or the {\it holomorphic
cotangent bundle} $T^\p M=\mathcal V(M)^\perp$, then $D$ is just
$\bar\partial_b$. If $M$ is a complex manifold and $E$ is an
analytic vector bundle, $D$ is $\bar
\partial$.

\section{The CR vector bundle structure on the holomorphic
cotangent bundle}

\subsection{Definition} Let now $M$ be an {\em integrable} CR manifold of CR
dimension $n = \dim \crb{M}_p$ and CR codimension $d=\dim M-2n$.
Integrability means that at each point $p\in M$ there exist a
neighbourhood $U$ of $p$ and  $N = n+d$ smooth CR functions $Z_j$
on $U$, $j=1,\dots ,N$ such that the differentials $dZ_1, \dots ,
dZ_N$ are linearly independent on $U$. The functions
$Z_1,\dots,Z_N$ are called a {\em family of basic solutions} in
$U$. One can show that $q\mapsto (Z_1(q),\ldots, Z_N(q))$ embeds
$M$, locally near $p\in M$, as a generic submanifold in $\C^N$.
Conversely, any embeddable (in this sense) CR manifold is
integrable.

We recall the definition of the holomorphic cotangent bundle $T^\p
M$: At each $p\in M$, $T^\p_p M = \crb{M}_p^\perp$. Since $\dim \C
T_p M = 2n + d$ and $\dim \crb{M}_p = n$, $\dim T^\p_p M = N$.
Hence, on a neighbourhood $U$ of a fixed point $p$ as above, the
$1$-forms $dZ_1, \dots, dZ_N$ form a basis of the sections of
$T^\p M$ over $U$. We want to define a CR atlas on $T'M$ such that
the bases obtained in this way are CR bases of sections. Assume
that we have another family of basic solutions $\td Z_1, \dots ,
\td Z_N$ on the open set $V\in M$. Then $d\td Z_j = \sum_k
a_j^{\;\,k} \, d Z_k$ on $U\cap V$, and exterior differentiation
yields
\[ \sum_k d a_j^{\;\,k} \wedge d Z_k = 0
\]
in the overlap $U\cap V$. By Cartan's Lemma, we conclude that $d
a_j^{\;\,k} = \sum C^l_k\, dZ_l$, and hence, ${\bar L} a_j^{\;\,k}
= \br{d a_j^{\;\,k}, {\bar L}} = 0$ for all CR vector fields
${\bar L}$, so that the $a_j^{\;\,k}$ are CR functions on $U\cap
V$. Hence, if we define a CR chart for $T^\p M$ over $U$ by
$\phi(\om_p) = (p, \ld^1,\dots,\ld^N)$ where $\om_p \in T^\p_p M$
is written as $\om_p = \sum_j \ld^j \,(dZ_j)_p$, then, by the
computation above, the collection of all such charts defines a CR
atlas on $T^\p M$. We equip $T'M$ with the canonical CR structure,
relative to this CR atlas, given by Lemma \ref{L:CRstruc}.

By the same argument as above using Cartan's lemma, we see that
the following lemma holds.

\begin{lem}\label{L:diffCR}
Let $g$ be a smooth function on an open set $U\subset M$. Then
$g$ is CR if and only if $dg$ is a CR section of $T^\p  M$ over
$U$.
\end{lem}

Sections of $T^\p M$ are called {\em holomorphic forms} (on $M$).
Let now $\om\in\Gm(M,T^\p M)$, and let ${\bar L}$ be a CR vector
field on $M$. The following useful formula holds:
\begin{equation}\label{E:DonTprime}
\CRder{{\bar L}} \om = {\bar L} \lrcorner d\om.
\end{equation}
To verify this identity, let $X\in\Gm (M,\C TM)$ be a tangent
vector field, and $\om = \sum_j a^j\, dZ_j$ locally. Then
\begin{align*}
\br{ {\bar L} \lrcorner d\om , X } &= \br{d\om , {\bar L}\wedge X} \\
&= \sum_j \br{ da^j \wedge dZ_j , {\bar L}\wedge X} \\
&= \sum_j (da^j({\bar L})\, dZ_j(X) - da^j (X)\, dZ_j ({\bar L})) \\
&= \sum_j ({\bar L}a^j)\, dZ_j (X) \\
&= \br{\CRder{{\bar L}}\om,X}.
\end{align*}

\subsection{Nondegeneracy conditions}\label{SS: conditions}
Recall that the {\em characteristic bundle} $T^0 M$ is defined by
\[ T^0 M = \left(\crb{M} \oplus \overline{\crb{M}}\right)^\perp.
\]
A {\em characteristic form} is a real nonvanishing section of $T^0
M$, and we write $\xi^0 (U)$ for the space of characteristic
forms on $U\subset M$. We usually denote a characteristic form by
$\gth$.

By the discussion above, an integrable CR manifold $M$ is finitely
nondegenerate (at $p$) as introduced by Baouendi, Huang and
Rothschild \cite{BHR1} if and only if for some $k$,
\[ \spanc \{ (\CRder{{\bar L}_1} \dots \CRder{{\bar L}_l} \gth) (p)
\colon {\bar L}_j \in \Gm(U,\crb{M}), \gth \in \xi^0 (U), l\leq k
\} = T^\p_p M,
\]
where $U$ is a suitable small neighbourhood of $p$. This follows
from the characterization of finite nondegeneracy found in e.g.
\cite{BERbook} and equation \eqref{E:DonTprime}.

We will use the notion of the pullback bundle introduced before in
order to give an intrinsic characterization of the nondegeneracy
of a CR map $f:M\to \hat M$ as introduced in \cite{L1}. Let $E =
f^\sharp T^\p \hat M$. $f$ is nondegenerate at $p_0$ if and only
if for some $k$,
\[
\spanc \{ \CRder{{\bar L}_1}\dots \CRder{{\bar L}_k} f^\sharp \gth
(p_0) \colon {\bar L}_j \in \Gm(U,\crb{M}), \gth \in \hat \xi
(\hat U), j\leq k \} = E_{p_0},
\]
where again $U$ is some neighbourhood of $p_0$, contained in
$f^{-1}(\hat U)$. More generally, we can define subspaces $E^k_p
\subset E_p$ (for $p\in U$, say) by
\[
E^k_p = \spanc \{ \CRder{{\bar L}_1}\dots \CRder{{\bar L}_l} \gth
(p) \colon {\bar L}_j \in \Gm(U,\crb{M}), \gth \in \xi^\p (U),
l\leq k \}
\]
Then $E^0_p \subset E^1_p \subset \dots \subset E^k_p \subset
\dots$, and for some smallest integer $k_0=k_0(p)$, $E^k_p =
E^{k_0}_p$ for $k \geq k_0 $. We call $s=s(p) = \dim E_{p} - \dim
E^{k_0(p)}_p$ the degeneracy of $f$ at $p$ and say that $f$ is
$(k_0,s)$-degenerate at $p$; we say that $f$ is of constant
degeneracy $s$ at $p_0$ if $s(p)\equiv d$ is constant in a
neighbourhood of $p_0$, and in this case we also say that $f$ is
constantly $(k_0, s)$-degenerate (note that the constancy only
applies to $s$). This characterization is equivalent to the
definition given in \cite{L2},  by a computation in coordinates
which we leave to the reader.

Note that $p\mapsto s(p)$ is upper semicontinuous. Hence, we have
the following:

\begin{lem} \label{L:global}
Let $f\colon M\to \hat M$ be a smooth CR map. Then there exists an
open dense subset $N\subset M$ such that $f$ is of constant
degeneracy at every point of $N$.
\end{lem}

Let us actually observe some more. It is proved in \cite{L2} that
the degeneracy $s$ of a transversal map is bounded by $\hat N -
N$. Furthermore it follows from that paper (and is easy to
prove), that $\dim E^1_p = N$ for any point $p$, if $\hat M$ is
Levi nondegenerate and $f$ is an embedding. Now note that if on
any open subset $U$ we had $E^{k+1}_p \subset E^k_p$ for $p\in
U$, it follows that $E_p = E^k_p$. It follows that if $f$ is of
constant degeneracy $s$ on $U$, then generically in $U$, $k_0
\leq \hat N - N + 1 -s \leq \hat N - N +1$. We state this as a
Lemma:

\begin{lem} \label{L:bounds}
Let $f\colon M\to \hat M$ be a smooth transversal CR embedding
between Levi nondegenerate hypersurfaces $M$ and $\hat M$. Then
there exists an open dense subset $N\subset M$ such that $f$ is
of constant degeneracy $(k_0,s)$ at every point of $N$ for some
$s\leq \hat N - N$ and $k_0\leq \hat N - N +1 - s \leq \hat N - N
+1$.
\end{lem}

\subsection{Some computations} From now on, unless stated otherwise, we assume
that
$M$ and $\hat M$ are of hypersurface type and integrable (so that
the holomorphic cotangent bundles are CR vector bundles), and that
$f:M\to \hat M$ is CR. Let us now fix a point $p_0$ (which we
denote by $0$ from now on) and a local basis ${L}_{\bar 1},\dots
{L}_{\bar n}$ of the CR vector fields on $M$ near $0$. We also
choose a characteristic form $\gth$ for $M$ near 0 and a real
vector field $T$ with $\br{\gth, T} =1$. We denote by $\gth^a$ the
dual form (relative to the basis of vector fields $T,L_1,\ldots,
L_n,L_{\bar 1},\ldots, L_{\bar n}$) to the anti CR vector field
$L_a$. We write $\hat 0 = f(0)$ and, as before, $\hat \gth$
denotes a (local) characteristic form on $\hat M$ near $\hat 0$,
and $\hat T$ a real vector field with $\br{\hat \gth, \hat T} =1$;
we extend $\hat\gth$ to a basis of $T^\p \hat M$ by $\hat n$ forms
$\hat \gth^1,\dots ,\hat \gth^{\hat n}$ which are the dual forms
of a local basis $\hat L_1,\dots \hat L_{\hat n}$ of the anti CR
vector fields on $\hat M$. We work on a small neighbourhood of $0$
which we may shrink, if necessary, without further mentioning it.

Since $f$ is CR, we can write
\[ f_* L_a = \gm_a^A \hat L_A,\quad  f_* L_{\bar a} = \overline{\gm_a^A} \hat
L,
\quad f_* T = \eta^A \hat L_A + \overline{\eta^A} \hat L_{\bar A}
+ \xi \hat T,\] or equivalently, \[ f^* \hat \gth^A = \gm^A_a
\gth^a + \eta^A  \gth, \quad f^* \hat \gth^{\bar A} =
\overline{\gm^A_a} \gth^{\bar a} + \overline{\eta^A}  \gth,\quad
f^* \hat \gth = \xi \gth.\] Here and from now on we use the
summation convention; small indices $a,b$ etc.\ range from $1$ to
$n$, capital indices $A,B$ etc.\ from $1$ to $\hat n$. We observe
here that $f$ is transversal at $p\in M$, as defined in the
introduction, if and only if $\xi(p)\neq 0$.

Let us write $D_{\bar a} = D_{L_{\bar a}}$, $f^\sharp \hat \gth =
\tau$, and $f^\sharp \hat \gth^A = \tau^A$; thus,
$\tau,\tau^1,\ldots,\tau^n$ forms a local basis (not necessarily
CR) for the CR vector bundle $E=f^\sharp T'\hat M$. We define
functions $g_{\bar a_1\dots \bar a_k B}$, $g_{\bar a_1 \dots \bar
a_k}$ on $M$ near $0$ by
\[D_{\bar a_k} D_{\bar a_{k-1}} \dots D_{\bar a_1} \tau =
g_{\bar a_1\dots \bar a_k B} \tau^B + g_{\bar a_1 \dots \bar a_k}
\tau.\]

Note that for elements $\sg_p$ of $E_p$, there is a natural
pairing $\br{\sg_p, X_{f(p)}}$ with elements $X_{f(p)}$ of $\C
T_{f(p)} \hat M$.  With this notation, we can write
\[
g_{\bar a_1\dots \bar a_k B} = \br{D_{\bar a_k} D_{\bar a_{k-1}}
\dots D_{\bar a_1} \tau , \hat L_B}, \quad
g_{\bar a_1 \dots \bar
a_k } = \br{ D_{\bar a_k} D_{\bar a_{k-1}} \dots D_{\bar a_1} \tau
, \hat T}.\]

We introduce similar notation for sections of $T^\p \hat M$ and
$T^\p M$:
\[ D_{\bar a_k} \dots D_{\bar a_1}  \gth = h_{\bar a_1 \dots \bar
a_k b} \gth^b + h_{\bar a_1\dots \bar a_k} \gth,\]
\[ D_{\bar A_k} \dots D_{\bar A_1}  \hat \gth = \hat h_{\bar A_1 \dots \bar
A_k B} \hat \gth^b + \hat h_{\bar A_1\dots \bar A_k} \hat \gth.\]
With $\br{\cdot,\cdot}$ denoting the canonical pairing between
vectors and covectors,
\[h_{\bar a_1\dots \bar a_k b} = \br{D_{\bar a_k} D_{\bar a_{k-1}}
\dots D_{\bar a_1} \gth ,  L_b}, \quad h_{\bar a_1 \dots \bar a_k
} = \br{ D_{\bar a_k} D_{\bar a_{k-1}} \dots D_{\bar a_1} \gth ,
 T},
\]
\[\hat h_{\bar A_1\dots \bar A_k B} = \br{D_{\bar A_k} D_{\bar A_{k-1}}
\dots D_{\bar A_1} \hat \gth , \hat L_B}, \quad \hat h_{\bar A_1
\dots \bar A_k } = \br{ D_{\bar A_k} D_{\bar A_{k-1}} \dots
D_{\bar A_1} \hat \tau , \hat T}.
\]

The $h_{\bar a b}$ relate to the Leviform of $M$ by
\begin{equation}\label{E:habasLeviform}
  h_{\bar a b} = \br{d\gth , L_{\bar a}\wedge L_b} = -\gth \left(
  [L_{\bar a}, L_b]
  \right),
\end{equation}
and similarly for $\hat M$.

We need to introduce (but just for a short while) one more
notation for our three bundles:
\[ D_{\bar a} \gth^b = R_{\bar a c}^b \gth^c + R_{\bar a}^b \gth,\]
\[D_{\bar a} \tau^B = {S}_{\bar a C}^B \tau^C + {S}_{\bar a}^B \tau,\]
\[D_{\bar A} \hat \gth^B = {\hat R}_{\bar A C}^B \hat \gth^C+ {\hat R}_{\bar
A}^B \hat \gth.\]

\begin{lem}\label{L:sections}
 Let $\sg$ be a section of $T^\p \hat M$. Then
\begin{equation}\label{E:transformationlaw1}
D_{\bar a} f^\sharp \sg = \overline{\gm_a^A} f^\sharp D_{\bar A}
\sg.
\end{equation}
In particular,
\begin{equation}\label{E:transformationlaw2}
  g_{\bar a B} = \overline{\gm_a^A} \hat h_{\bar A B},
  \quad g_{\bar a} = \overline{\gm_a^A} \hat h_{\bar
  A},
\end{equation}
and
\begin{equation}\label{E:transformaionlaw3}
  {S}_{\bar a C}^B = \overline{\gm_a^A} {\hat R}_{\bar A C}^B,
  \quad {S}_{\bar a}^C = \overline{\gm_a^A} {\hat R}_{\bar A}^a.
\end{equation}
\end{lem}
\begin{proof}
Let us consider a local CR basis $\om_1,\dots \om_{\hat N}$ for
$T^\p \hat M$. We write $\al_j = f^\sharp \om_j$ for the
corresponding local CR basis of $E$. Then locally $\sg = a^j
\om_j$ for smooth functions $a^j$, and we have
\begin{align*}
D_{\bar a} f^\sharp \sg &= D_{\bar a} f^\sharp (a^j \om_j)\\
&= D_{\bar a} (a^j\circ f) \al_j\\
&= (L_{\bar a} (a^j \circ f )) \al_j \\
&= ((f_* L_{\bar a}) a^j) \al_j \\
& = \overline{\gm_a^A} (L_{\bar A} a^j \circ f) \al_j \\
&= \overline{\gm_a^A} f^\sharp D_{\bar A} \sg.
\end{align*}

For the other assertions, we note that
\[ D_{\bar a} \tau = g_{\bar a B} \tau^B + g_{\bar a} \tau \] and
\[ D_{\bar A} \hat \gth = \hat h_{\bar A B} \hat \gth^B + h_{\bar
A} \hat \gth. \]%
Since
\[ \overline{\gm_a^A} f^\sharp (D_{\bar A} \hat \gth) =
\overline{\gm_a^A} \hat h_{\bar A B} \tau^B + \overline{\gm_a^A}
\hat h_{\bar A} \tau,\] comparing coefficients we get
\eqref{E:transformationlaw2}. \eqref{E:transformaionlaw3} follows
by the same argument.
\end{proof}

The following lemma is an immediate consequence of
\eqref{E:productrule} and the definitions above:
\begin{lem}\label{L:derivativrule} With the notations introduced above, the
following hold:
\begin{equation}\label{E:derivativrule1}
  g_{\bar a_1 \dots \bar a_k \bar a B} = L_{\bar a}g_{\bar a_1 \dots \bar a_k
  B} + g_{\bar a_1 \dots \bar a_k  C} S_{\bar a B}^C + g_{\bar a_1 \dots \bar
a_k }g_{ \bar a B}
\end{equation}
\begin{equation}\label{E:derivativrule2}
  {\hat h}_{\bar A_1 \dots \bar A_k \bar A B} = L_{\bar A} {\hat h}_{\bar A_1
\dots \bar A_k
  B} + {\hat h}_{\bar A_1 \dots \bar A_k  C} {\hat R}_{\bar A B}^C + {\hat
h}_{\bar A_1 \dots \bar A_k }{\hat h}_{ \bar A B}
\end{equation}
\begin{equation}\label{E:derivativrule3}
  h_{\bar a_1 \dots \bar a_k \bar a b} = L_{\bar a}h_{\bar a_1 \dots \bar a_k
  b} + h_{\bar a_1 \dots \bar a_k  c} R_{\bar a b}^c + h_{\bar a_1 \dots \bar
a_k }h_{ \bar a
  b}.
\end{equation}
\end{lem}
We also record the following two lemmas for later use. We write
$I= (i_1\dots i_l)$, where $1\leq i_j \leq n$, and we use the
notation $L^I = L_{i_1} \dots L_{i_l}$; for such an $I$, we write
$|I| = l$.
\begin{lem}\label{L:degeneracyincoords} For each $k\in\N$, the two subspaces
of $\C^{\hat n}$ defined by
\[ \spanc \{(g_{\bar a_1 \dots \bar a_l A} (p))_A \colon l\leq k
\}
\]
\[\spanc \{(L^{\bar I} g_{\bar a A} (p))_A \colon |I|\leq k
\}
\]
coincide. Consequently, the following are equivalent definitions
for the degeneracy $d(p)$ of $f$ at $p$ defined in section
\ref{SS: conditions}:
\begin{enumerate}[(i)]
\item $d(p) = \hat n - \max_k \dim \spanc \{(g_{\bar a_1 \dots \bar a_k A}
(p))_{A} \colon l\leq k\}$;
\item $d(p) = \hat n - \max_k \dim \spanc \{ (L^{\bar I} g_{\bar a A} (p))_A
\colon |I|\leq k\}$;
\end{enumerate}
 and in each case, $k_0 (p)$ is the minimal integer
$k$ for which the maximum on the right hand side is attained.
\end{lem}
The proof is by induction using \eqref{E:derivativrule1}.
\begin{lem}\label{L:gisequalrational}
For each $a_1,\dots, a_k$, $1\leq a_l\leq n$, and $B$, $1\leq
B\leq \hat n$ there exist  functions $u_{\bar a_1\dots \bar a_k
B}$ and $u_{\bar a_1\dots \bar a_k}$ polynomial in their arguments
with coefficients smooth on $M\times\hat M$  such that
\begin{equation}\label{E:gisequalrational1}
    g_{\bar a_1\dots \bar a_k B} = u_{\bar a_1\dots \bar a_k B}
    \left( \overline{L^I \gm_a^A} \right),\quad  {|I| \leq k -1}
\end{equation}
\begin{equation}\label{E:gisequalrational2}
    g_{\bar a_1\dots \bar a_k } = u_{\bar a_1\dots \bar a_k}
    \left( \overline{L^I \gm_a^A}  \right), \quad {|I| \leq k -1}.
\end{equation}
Furthermore,  $u_{\bar a_1\dots \bar a_k A}$ and $u_{\bar a_1\dots
\bar a_k}$ depend only on $M$ and $\hat M$ (and not on the
mapping $f$).
\end{lem}
The equalities \eqref{E:gisequalrational1} and
\eqref{E:gisequalrational2} mean the following: if the left hand
side is evaluated at $p$, it equals the right hand side with the
coefficients evaluated at $(p,f(p))$ and the arguments evaluated
at $p$. We will encounter many such equalities, and we will
always use this time-saving notation.

For the proof of Lemma~\ref{L:gisequalrational}, we recall that
for $k=1$, $g_{\bar a B} = \overline{\gm_a^A} \hat h_{\bar A B}$
and $g_{\bar a} = \overline{\gm_a^A} \hat h_{\bar A}$ by
Lemma~\ref{L:sections}. An induction and
Lemma~\ref{L:derivativrule}  finishes the proof of the
Lemma~\ref{L:gisequalrational}; the details are left to the
reader.

\subsection{The pullback operation}\label{SS:pullback}
In this section, we discuss
the relation between our vector bundle $E=f^\sharp T'\hat M$ and
the vector bundle $T^\p M$ over $M$. Recall that we have a pairing
$\br{\sg, X}$ for sections of $E$ and tangent vectors $X$ at
$f(p)$. We define a map $\phi$ by
\[
\br{\phi \sg, X} = \br{\sg,f_* X}.
\]
Since $f$ is CR, $\phi$ maps sections of $E$ to sections of $T^\p
M$. Note that $\phi: \Gm(U,E) \to \Gm(U,T^\p M)$ is linear over
$C^\infty (U)$, since
\[\br{\phi a \sg, X} = (a \sg)(f_*X) = a \sg(f_*X) = a \br{\phi
\sg, X}.
\]
\begin{lem}\label{L:commutativity} For each $j$, $D_{\bar j} \phi = \phi
D_{\bar
j}$. Furthermore, $\phi \circ f^\sharp = f^*$.
\end{lem}
\begin{proof}
 Since the derivative operation is local,
it is actually enough to work on a trivializing open set $U=
f^{-1} (V)$, where  have a family of basic solutions $\hat Z_j$
in $V$. A basis of the CR sections of $E$ is then given by $\sg_j
= f^\sharp d \hat Z_j$. Let $\tau$ be any section of $E$ over $U$.
Then $\tau = a^j \sg_j$, and $\phi \tau = a^j \phi(\sg_j)$. By the
definition of $D_{\bar j}$ we now see it is enough to show that
if $\sg$ is a CR section of $E$, then $\phi\sg$ is a CR section
of $T^\p M$. We compute:
\begin{align*} \br{ (\phi f^\sharp d\hat Z_j)_p , X_p}
&= \br{ (f^\sharp d\hat Z_j)_p , f_* X_p} \\
&= \br{ (d\hat Z_j)_{f(p)}, f_* X_p}  \\
&= \br{ f^* (d \hat Z_j)_p, X_p} \\
&= \br{ d( \hat Z_j\circ f)_p, X_p},
\end{align*}
so that actually $\phi f^\sharp d\hat Z_j = d (\hat Z_j\circ f)$,
which is CR by Lemma \ref{L:diffCR}. The second claim follows from
the same calculation.
\end{proof}
By an appropriate choice of basis $\gth^a$, $\gth$, and similarly
on $\hat M$, we may assume that $R^b_{\bar a c}, R^b_{\bar a},
\hat R^B_{\bar A C}, R^B_{\bar A}$ are all $0$ (cf.\ \cite{E5},
Lemma 1.33). By Lemma \ref{L:sections}, it follows that $S^B_{\bar
a c}=S^B_{\bar a}=0$. Lemma \ref{L:commutativity} implies the
following relations:
\begin{equation}\label{E:levitransform}
  \xi h_{\bar a b} = \gm_a^B g_{\bar b B};
\end{equation}
\begin{equation}\label{E:xitransform}
 L_{\bar a} \xi + \xi h_{\bar a} = g_{\bar a B} \eta^B + g_{\bar a
 } \xi.
 \end{equation}
This follows comparing the coefficients of $D_{\bar a} \phi \tau$
and $\phi D_{\bar a } \tau$. We also need the following
equalities:
\begin{equation}\label{E:barredderivativeofgamma}
  L_{\bar a} \gm_{b}^B = - \eta^B h_{\bar a b},
\end{equation}
\begin{equation}\label{E:barredderivativeofeta}
  L_{\bar a} \eta^B = - \eta^B h_{\bar a},
\end{equation}
\begin{equation}\label{E:Tderivativeofgamma}
  T\gm_a^A = L_a \eta^A + \eta^A h_a,
\end{equation}
which follow as in \cite{E5}. We observe that the identity
\eqref{E:levitransform} implies that if $f$ is transversal at $p$
(i.e.\ $\xi(p)\neq 0$) and $M$ is Levi nondegenerate at $p$ (i.e.
the matrix $(h_{\bar a b})(p)$ is invertible), then the rank of
the matrix $(\gamma^B_a)(p)$ must be maximal, i.e. $=n$. Thus, if
$f$ is transversal and $M$ is Levi nondegenerate, then $f$ is in
fact a local embedding.

\section{Jet reflection identities for nondegenerate mappings }
In this section we are going to establish the jet reflection
identities necessary to use the machinery developed in \cite{E5}.
Our first goal is the reflection identities for nondegenerate
mappings. Fix an integer $k_0$. We denote by $\indexset$ the set
of all strings $\stringA = \bar a_1 \dots \bar a_k$ of elements of
$\{\bar 1,\dots, \bar n\}$ with $k\leq k_0$, and we write
$\indexsettwo$ for the set of all $\hat n$-tuples of elements of
$\indexset$. For $\stringA = \bar a_1 \dots \bar a_k$ we let
$\ell(\stringA) = k$ denote the length of $\stringA$.

\begin{thm}\label{T:reflectionnondegenerate} Let $k_0\in \N$. There exist
functions
\begin{equation*}
\Da_{\dblind}\left( \overline{ L^J \gm_c^D }\right) , \quad
r_{a,\dblind}^{B} \left( \overline{ L^J \gm_c^D }, \overline{L^I
\xi}
  \right), \quad r^{B}_{\dblind}\left( \overline{ L^J \gm_c^D }, \overline{L^I
\xi}
  \right),
\end{equation*}
where  $|J|\leq k_0 - 1, |I| \leq k_0$, for $\dblind =
(\mathbf{A}_0,\dots,\mathbf{A}_n) \in\indexsettwo^{n+1}$, $a\in \{
1,\dots ,n\}$ and $B\in \{1,\dots,\hat n\}$ with the following
property: If $f:M\to \hat M$ is $k$-nondegenerate at $0$ with
$k\leq k_0$, then there exists $\dblind^0$ with $\Da_{\dblind^0}
\left( \overline{L^J \gm_c^D} (0) \right) \neq 0$, and for any
$\dblind \in\indexsettwo^{n+1}$ with $\Da_{\dblind} \left(
\overline{L^J \gm_c^D} (0) \right) \neq 0$ the following
equalities hold:
\begin{equation}\label{E:reflectionnondegenerate1}
  \gm_a^B = r_{a,\dblind}^B \left( \overline{ L^J \gm_c^D }, \overline{L^I
\xi}
  \right),
\end{equation}
\begin{equation}\label{E:reflectionnondegenerate2}
  \eta^B = r^B_\dblind \left( \overline{ L^J \gm_c^D }, \overline{L^I \xi}
  \right).
\end{equation}
The $\Da_\dblind$ are polynomials in $ \overline{L^J \gm_c^D}$
with complex coefficients; the functions $r_{a,\dblind}^B$ and
$r_\dblind^B$ are rational in $ \overline{L^J \gm_c^D}$,
polynomial in $\overline{L^I\xi}$, with coefficients smooth on
$M\times\hat M$. Furthermore, $\Da_\dblind$, $r_{a,\dblind}^B$,
and $r_\dblind^B$ only depend on $M$, $\hat M$, and $k_0$ (and
not on the mapping $f$).
\end{thm}

In order to define the functions $\Da_\dblind$ of the Theorem, we
need the following computational Lemma.
\begin{lem}\label{L:complemma} Let $f:M\to \hat M$ be a CR mapping.
Then for each $k\in \N$, all $a \in \{1,\dots, n\}$, and all
$\stringA = \bar a_1 \dots \bar a_k$ the following identities
hold:
\begin{equation}\label{E:complemma1}
  g_{\stringA B}\eta^B =
  \sum_{\ell(\stringA^\p) < \ell(\stringA) } g_{\stringA^\p B} \eta^B \,
v_{\stringA}^{\stringA^\p}
    \left( \overline{L^J \gm_c^D}\right) + w_{\stringA} \left( \overline{L^J
\gm_c^D},\overline{L^I \xi} \right);
\end{equation}
\begin{multline}\label{E:complemma2}
     g_{\stringA B} \gm_a^B =
  \sum_{\ell(\stringA^\p) < \ell(\stringA)}
  g_{\stringA^\p B} \gm_a^B \, p_{a \stringA}^{\stringA^\p}
    \left( \overline{L^J \gm_c^D}\right) + \\
  \sum_{\ell(\stringA^\p) < \ell(\stringA) } g_{\stringA^\p} \eta^B \, q_{a
\stringA}^{\stringA^\p}
    \left( \overline{L^J \gm_c^D}\right) + r_{a \stringA} \left( \overline{L^J
\gm_c^D},\overline{L^I \xi} \right);
\end{multline}
where $|J|\leq k-1$, $|I|\leq k$. The functions $v$, $w$, $p$,
$q$, and $r$ are rational in $\overline{L^I \gm_c^D}$, the
functions $w$ and $r$ are polynomial in $\overline{L^I \xi}$, with
coefficients smooth on $M\times\hat M$. Furthermore, $v$, $w$,
$p$, $q$, and $r$ only depend on $M$ and $\hat M$.
\end{lem}
\begin{proof} We start with \eqref{E:complemma1}. For $k=1$, by
\eqref{E:xitransform},
\[g_{\bar a B} \eta^B = L_{\bar a} \xi + \xi h_{\bar a} - \xi
g_{\bar a},
\]
and recalling \eqref{E:transformationlaw2}, $ g_{\bar a} =
\overline{\gm_a^A} h_{\stringA} $, we have that
\[ g_{\bar a B} \eta^B =
L_{\bar a} \xi + \xi h_{\bar a} - \xi
\overline{\gm_a^A}h_{\stringA}= w_{\bar a} \left(
\overline{\gm_c^D}, \overline{L^I \xi} \right),
\]
which is \eqref{E:complemma1} for $k = 1$. To get
\eqref{E:complemma1} for $k>1$, we use \eqref{E:complemma1} for
$k-1$ to express $g_{\bar a_1 \dots \bar a_{k-1} B} \eta^B$ and
apply $L_{\bar a_k}$, using Lemma~\ref{L:derivativrule} and
\eqref{E:barredderivativeofeta}.

The proof of \eqref{E:complemma2} is similar: We start with
\eqref{E:levitransform} for the case $k=1$, and for $k>1$, we use
induction, Lemma~\ref{L:derivativrule},
\eqref{E:barredderivativeofgamma}, and
\eqref{E:barredderivativeofeta}.
\end{proof}
\begin{proof}[Proof of Theorem~\ref{T:reflectionnondegenerate}]
Let $\fatA = (\stringA_1,\dots, \stringA_{\hat n})$. We rewrite
the $\hat n$ equations obtained from  \eqref{E:complemma1} by
replacing  $\stringA$ with $\stringA_1,\dots, \stringA_{\hat n}$
in the following form:
\begin{equation}\label{E:systemnondegenerateforeta}
 M_{\fatA} \eta =
w_{\fatA},
\end{equation}
where
\[\eta = \begin{pmatrix}
  \eta^1 \\
  \vdots \\
  \eta^{\hat n}
\end{pmatrix},\quad
 w = \begin{pmatrix}
  w_{\stringA_1} \\
  \vdots \\
  w_{\stringA_{\hat n}}
\end{pmatrix},
\]
and the $(i,B)$th entry of $M_{\fatA}$ is
\[g_{\stringA_i B} -
  \sum_{\ell(\stringA^\p) < \ell(\stringA_i)}
  g_{\stringA^\p B} \, v_{\stringA_i}^{\stringA^\p}.
\]
Similarly, for each $a\in\{1,\dots,n\}$ we rewrite
\eqref{E:complemma2} as
\begin{equation}\label{E:systemnondegenerateforgamma}
N_{a,\fatA} \gm_a = K_{a,\fatA} \eta + r_\fatA.
\end{equation}
By Lemma~\ref{L:gisequalrational} and Lemma~\ref{L:complemma} the
matrices $M_\fatA$ and $N_{a,\fatA}$ are polynomial in
$\overline{L^J \gm_c^D}$, with coefficients smooth on $M\times
\hat M$ and they only depend on $M$ and $\hat M$. We define
$\Da_{0,\fatA} = \det M_\fatA$ and $\Da_{j,\fatA} = \det
N_{j,\fatA}$, $j=1,\dots, n$, where we evaluate the coefficients
of $\overline{L^J \gm_c^D}$ at $(0,0)$, and for $\dblind =
(\fatA_0, \dots ,\fatA_n)$, we set $\Da_{\dblind}=\prod_{j=0}^n
\Da_{j,\fatA_j}$. By Cramer's rule, if $\Da_{\dblind}
\left(\overline{L^J \gm_c^D}(0) \right) \neq 0$, in a
neighbourhood of $0$ we can solve the equations
\eqref{E:systemnondegenerateforeta} with $\fatA = \fatA_0$ and
\eqref{E:systemnondegenerateforgamma} with $(a, \fatA) =
(a_1,\fatA_1),\dots,(a_n,\fatA_n)$ for $\eta^B$ and $\gm_a^B$ in
the required form.

We finish the proof by showing that if $f$ is $k$-nondegenerate
at $0$  with $k\leq k_0$, then there exists $\dblind^0$ with
$\Da_{\dblind^0} \left(\overline{L^J \gm_c^D}(0) \right) \neq 0$.
First recall that by Lemma~\ref{L:degeneracyincoords}, $f$ is
$k$-nondegenerate at $0$ if and only if
\begin{equation}\label{E:nondegenerateifandonlyif}
\spanc \{ (g_{\stringA B}(0))_B  \colon \ell(\stringA)\leq k \} =
\Cnh.
\end{equation}
Consider a system of vectors of the form
\begin{equation}\label{E:definitionofthevs}
  v_{\stringA} =  \left(g_{\stringA B} (0) - \sum_{\ell(\stringA^\p) <
\ell(\stringA )}
c_{\stringA B}^{\stringA^\p} g_{\stringA^\p B} (0)\right)_B, \quad
\ell(\stringA) \leq k.
\end{equation}
We claim that if $f$ is $k$-nondegenerate, then $\spanc \{
v_{\stringA} \colon \ell(\stringA)\leq k \} = \Cnh$, which implies
that $\Da_{\dblind^0} \left(\overline{L^J \gm_c^D}(0) \right) \neq
0$ for some $\dblind^0$. To prove the last claim, we let $V_l =
\spanc \{(g_{\stringA B})_B (0) \colon \ell(\stringA)\leq l \}$
and $W_l = \spanc \{ v_{\stringA} \colon \ell(\stringA)\leq l \}$.
By \eqref{E:definitionofthevs}, $V_1 = W_1$. But then, again by
\eqref{E:definitionofthevs}, $V_2 = W_2$. By induction, we see
that $W_k = V_k = \Cnh$.
\end{proof}

Theorem \ref{T:reflectionnondegenerate} is completely analogous to
Theorem 2.4 in \cite{E5}. The proof of Theorem
\ref{T:nondegenerate} can now easily be completed using the
arguments in \cite{E5} and a result from \cite{L2}.

\begin{proof}[Proof of Theorem~\ref{T:nondegenerate}]
By repeating the arguments in \cite{E5}, based on the identities
in Theorem~\ref{T:reflectionnondegenerate} instead of those in
Theorem 2.4 in \cite{E5}, we obtain a complete system of some
order $l_0$ (the exact value of which is not relevant for this
argument) for $f$ and $g$. That is, there exists a neighbourhood
$U$ of $p_0$ and a smooth function $\phi:U\times J^{l_0}(U,\hat M)
\to J^{l_0+1}(U,\hat M)$ such that
\[j^{l_0 + 1}_{x} f = \phi (x,j^{l_0}_x f) \text{ and }
 j^{l_0 + 1}_{x} g = \phi (x,j^{l_0}_x g), \quad x\in U.\] By
Theorem 2 in \cite{L2}, $j^l_{p_0} f = j^l_{p_0} g$ for all $l$
and, in particular then, $j^{l_0}_{p_0} f = j^{l_0}_{p_0} g$. The
conclusion of Theorem~\ref{T:nondegenerate} follows from the
uniqueness of solutions of complete systems (see e.g. \cite{BCG}).
\end{proof}
\section{Jet reflection identities for Levi-nondegenerate hypersurfaces}
In this section, we consider the case $\hat N = N + 1$, where $M$
is Levi-nondegenerate at $0$, and $\hat M$ is Levi-nondegenerate
at $\hat 0$. We are also assuming that $f$ is transversal, which
in our notation is equivalent to $\xi (0) \neq 0$.

For the proof of Theorem~\ref{T:codimension1} we first note that
by Lemma 20 in \cite{L2}, the degeneracy of $f$ is either $0$ or
$1$. In the case where the degeneracy of $f$ is $0$,
Theorem~\ref{T:codimension1} follows from
Theorem~\ref{T:nondegenerate}, so that we only have to deal with
the case where $f$ is constantly $1$-degenerate. Using
\eqref{E:levitransform} it is easy to see that with our
assumptions, the vectors  $g_{a B} (0)$, $ a=1,\dots, n$ are
linearly independent, as are the vectors $ \gm_c^D (0)$,  $
c=1,\dots, n$. We conclude that $f$ is constantly
$(1,1)$-degenerate. Our goal in this section is to derive the
following jet reflection identity for such mappings. The unique
determination theorem will be proved in
section~\ref{S:completesystem} by producing a complete system
using this jet reflection identity.
\begin{thm}\label{T:reflectionlevinondegenerate}
Let $M\subset\CN$ and $\hat M\subset\CNplus$ be hypersurfaces
which are Levi-nondegenerate at the points $p_0$ and $\hat p_0$,
$f\colon M\to \hat M$ a smooth CR map which is constantly
$(1,1)$-degenerate and transversal at $p_0$. Then for any $I$ with
$|I| = k$ there exist functions $r^{I,A}_b$ and $s^{I,A}$ for
$a,b\in \{1,\dots,n\}$ and $A\in \{1,\dots,n+1\}$ such that
\begin{align}\label{E:reflectionlevinondegenerategamma}
  L^I \gm_b^A &= r_{b}^{I,A} \left( L^M \gm_c^D, L^M \eta^D,
  \overline{\gm_c^D}, \overline{\eta^D} \right),\\
\label{E:reflectionlevinondegenerateeta}
  L^I \eta^A &= s^{I,A} \left( L^M  \gm_c^D,L^M \eta^D,
  \overline{L_d\gm_c^D}, \overline{L_d\eta^D},\overline{T\gm_c^D},
  \overline{T\eta^D} \right)
\end{align}
Here $|M| \leq k-1$. The $r^{I,A}_b$ and $s^{I,A}$ are rational in
their arguments with coefficients which are smooth on $M\times
\hat M$, and they only depend on $M$ and $\hat M$.
\end{thm}

\begin{proof} Since $\hat M$ is Levi-nondegenerate, we can choose the basis
of CR vector fields tangent to $\hat M$ such that
\begin{equation}\label{E:conventionlevinondegenerate}\hat h_{\bar A B} (0)
= \eps_{\bar A B} = \begin{cases} \pm i, & A = B \\
0, & A\neq B \end{cases}.\end{equation}

Since $f$ is constantly $(1,1)$-degenerate, for any ${\bar a}$,
$\bar b$ the determinant of the matrix
\begin{equation}\label{E:bigmatrix}
\begin{pmatrix}
  g_{\bar 1 \, 1} & \dots & g_{\bar 1 \, n+1} \\
  \vdots &  & \vdots \\
  g_{\bar n \, 1} & \dots & g_{\bar n \, n+1} \\
  L_{\bar a} g_{\bar b \, 1} & \dots  & L_{\bar a} g_{\bar b \, n+1}
\end{pmatrix}
\end{equation}
vanishes. Let us write $(-1)^{B+1} \Da^B$ for the determinant of
the matrix obtained from the matrix \eqref{E:bigmatrix} by
dropping the last row and the $B$th column. We fix $a$ and $b$
for the moment. By developing the determinant of
\eqref{E:bigmatrix} along the last row,
\begin{equation}\label{E:detdeveloped}
  \Da^B L_{\bar a} g_{\bar b B} = 0.
\end{equation}
We compute
\begin{equation}\label{E:Lbaragisequal}
  L_{\bar a} g_{\bar b B} = L_{\bar a}\left( \overline{\gm_b^A}
  {\hat  h}_{\bar A B}\right) =
   \overline{L_{ a}\gm_b^A} {\hat  h}_{\bar A B} +
     \overline{\gm_b^A \gm_a^D} {\hat L}_{\bar D}
     {\hat h}_{\bar A B}.
\end{equation}
Plugging \eqref{E:Lbaragisequal} into \eqref{E:detdeveloped} and
using Lemma~\ref{L:gisequalrational} we conclude that
\begin{equation}\label{E:equationforbarlgamma}
 \Da^B \overline{L_{ a}\gm_b^A} {\hat  h}_{\bar A B} = r\left(
 \overline{ \gm_c^D } \right),
\end{equation}
where $r$ is rational in $\overline{\gm_c^D}$ with coefficients
smooth on $M\times \hat M$ and $r$ is determined only by $M$ and
$\hat M$. We take the complex conjugate of this equation to obtain
\begin{equation}\label{E:equationforlgamma}
 \Da^{\bar B} {L_{ a}\gm_b^A} {\hat  h}_{ A \bar B} = {\bar r}\left(
  \gm_c^D  \right),
\end{equation}

Consider now the equations \eqref{E:levitransform} which we
rewrite as
\[\xi h_{\bar d b} = \gm_b^A g_{\bar d A}.
\]
We apply $L_a$ to this equation and conclude, using also equations
\eqref{E:levitransform} and \eqref{E:xitransform},
\begin{equation}\label{E:lalevitransform}
  \left( L_a \gm_b^A \right) g_{\bar d A} = \td p_{\bar d} \left( \gm_c^D,
  \overline{\gm_c^D}, \xi, L_a \xi \right) = p_{\bar d} \left(
  \gm_c^D,\eta^D
  \overline{\gm_c^D}, \overline{\eta^D} \right),
\end{equation}
where $p_{\bar d}$ is rational in its arguments with coefficients
smooth on $M\times \hat M$ and $p_{\bar d}$ is determined only by
$M$ and $\hat M$. We claim that the matrix
\begin{equation}\label{E:matrixforsystem}
  \begin{pmatrix}
    g_{\bar 1\, 1}(0)  & \dots  & g_{\bar 1 \, n+1}(0) \\
    \vdots &  & \vdots \\
    g_{\bar n \, 1}(0) & \dots & g_{\bar n \, n+1} (0) \\
    \Da^{\bar 1}(0) \eps_{\bar 1 1} & \dots & \Da^{\overline{n+1}} (0)
  \eps_{\overline{n+1}n+1}
  \end{pmatrix}
\end{equation}
is invertible, so that the system of equations
\begin{align}\label{E:systemforlevinondegenerate}
  X^A g_{\bar d A} &= p_{\bar d}, \quad d\in\{1,\dots ,n\},\\
  X^A \Da^{\bar B} {\hat h}_{ A \bar B} &= \bar
r,\label{E:systemforlevinondegenerate2}
\end{align}
can be solved uniquely for the $X^A$ (in a neighbourhood of $0$).
Developing the determinant $\Da$ for the matrix
\eqref{E:matrixforsystem} along the last row, we conclude that
\begin{equation}\label{E:deltaisequal}
  \Da = \eps_{\bar A A} \Da^{\bar A} \Da^A.
\end{equation}
We are now going to show that this last expression is not equal
to $0$. Consider the matrices
\begin{equation}\label{E:matricesGandGm}
  G= \begin{pmatrix}
    g_{\bar 1\, 1}(0)  & \dots  & g_{\bar 1 \, n+1}(0) \\
    \vdots &  & \vdots \\
    g_{\bar n \, 1}(0) & \dots & g_{\bar n \, n+1} (0)
  \end{pmatrix},
  \quad
  \Gm = \begin{pmatrix}
    \gm_{\bar 1}^{ 1}(0)  & \dots  & \gm_{\bar n}^{1}(0) \\
    \vdots &  & \vdots \\
    \gm_{\bar 1}^{ n+1}(0) & \dots & \gm_{\bar n}^{ n+1} (0),
  \end{pmatrix}
\end{equation}
and the $n\times n$ matrix
\begin{equation}\label{E:matrixH}
  H=\begin{pmatrix}
    h_{\bar 1\, 1}(0)  & \dots  & h_{\bar 1 \, n}(0) \\
    \vdots &  & \vdots \\
    h_{\bar n \, 1}(0) & \dots & h_{\bar n \, n} (0) \\
  \end{pmatrix}.
\end{equation}
In matrix notation, \eqref{E:levitransform} then reads as $\xi H
= G \Gm$. Let us write $G_A$ for the matrix $G$ with the $A$th
column dropped, and $\Gm^A$ for the matrix $\Gm$ with the $A$th
row dropped. The Cauchy-Binet Theorem from linear algebra states
that $\det G \Gm^t = \det G_A\det \Gm^A$. By our convention
\eqref{E:conventionlevinondegenerate} and
\eqref{E:transformationlaw2} it follows that
\[ g_{\bar a B} (0) = \overline{\gm_a^B}(0) \eps_{\bar B B}.
\]
Writing $\eps = \prod_B \eps_{\bar B B}$ we obtain from the
equation above
\[ \Da^A = \frac{\eps}{\eps_{\bar A A}} \det \overline{\Gm^A}.
\]
Recalling that $\Da^A = \det G_A$, we see that the Cauchy-Binet
Theorem implies that
\begin{align*}\label{E:deltaisequal2}
  \Da &= \eps_{\bar A A} \Da^{\bar A} \Da^A \\
  &= \eps \det {\overline{G_A}} \det \overline{\Gm^A} \\
  &= \eps \xi(0)^n \det H \neq 0,
\end{align*}
and our last claim is proved. Using that the matrix
\eqref{E:matrixforsystem} is invertible to solve
\eqref{E:equationforlgamma} and \eqref{E:lalevitransform} for $
L_a \gm^A_b$, we obtain (by repeatedly differentiating the
resulting equations and using the identities
\eqref{E:barredderivativeofgamma},
\eqref{E:barredderivativeofeta}, and \eqref{E:Tderivativeofgamma})
the first part of Theorem \ref{T:reflectionlevinondegenerate}.

We also need reflection identities
\eqref{E:reflectionlevinondegenerateeta} for $L_a \eta^B$. We will
deduce those from the equation
\eqref{E:reflectionlevinondegenerategamma} we just proved. We
first compute
\[ L_{\bar c} L_b \gm_e^B = L_b L_{\bar c} \gm_e^B - h_{\bar c b}
T\gm_e^B,
\]
which follows from \eqref{E:habasLeviform} and the fact that
$R^b_{\bar a c}=0$ (cf. \cite{E5}, Lemma 1.33). By
\eqref{E:barredderivativeofgamma} and
\eqref{E:Tderivativeofgamma}, we conclude that
\begin{equation}\label{E:new} L_{\bar c} L_b \gm_e^B = - L_b
\left( \eta^B h_{\bar c e} \right) - h_{\bar c b}\left( L_e \eta^B
+ \eta^B h_e \right)= -( h_{\bar c e} L_b \eta^B+h_{\bar c b} L_e
\eta^B) + r(\eta^B) ,\end{equation} where $r(\eta^B)$ is a
(linear) polynomial in $\eta^B$ with coefficients which are smooth
on $M$ (and only depend on $M$). The identity
\eqref{E:reflectionlevinondegenerateeta} with $k=1$ now follows by
solving for $L_b\eta^B$ in \eqref{E:new} (e.g.\ set $b=c=e$ and
recall that $h_{\bar b b}\neq 0$) and then applying $L_{\bar c}$
to \eqref{E:reflectionlevinondegenerategamma} with $k=1$, and
finally using \eqref{E:barredderivativeofgamma} and
\eqref{E:barredderivativeofeta} to substitute for $L_{\bar
c}\gamma^B_e$ and $L_{\bar c}\eta^B$. The general case follows by
differentiating. The details are left to the reader.
\end{proof}

\section{Strictly pseudoconvex targets}
Our goal in this section is to derive the following reflection
identity; in section~\ref{S:completesystem} we will use it to
derive a complete system for the mapping $f$.

\begin{thm}\label{T:reflectionstrpscx} Let $M\subset\CN$, $\hat M \subset\CNh$
be smooth
hypersurfaces, $p_0 \in M$, $\hat p_0 \in \hat M$, with $\hat M$
strictly pseudoconvex at $\hat p_0$. Also assume that $f$ is a
smooth CR mapping from $M$ to $\hat M$ with $f(p_0) = \hat p_0$
which is transversal and constantly $(k_0,s)$-degenerate at $p_0$.
Then for any $I$ with $|I| = k$, $a\in\{1,\dots,n\}$,
$B\in\{1,\dots,\hat n\}$ there exist functions $r_a^{I,B}$ and
$s^{I,B}$ which are rational in their arguments such that
\begin{align}\label{E:reflectionpseudoconvexgamma}
 L^I \gm_a^B &= r_a^{I,B} \left( L^M T^l \gm_c^D,L^M T^l \eta^D,
 \overline{ L^J T^m \gm_c^D}, \overline{ L^J T^m \eta^D}
 \right),\\ \label{E:reflectionpseudoconvexeta}
 L^I \eta^B &= s^{I,B} \left( L^M T^l \gm_c^D, L^M T^l \eta^D,
 \overline{ L^J T^m \gm_c^D}, \overline{ L^J T^m \eta^D}
 \right).
\end{align}
Here $l \in \{ 0,1 \}$, $|M| + l \leq \max (k_0,k-1)$, $|J| + m
\leq k_0+1$, $m\leq \min(k,k_0)$. Furthermore, $r_a^{I,B}$ and
$s^{I,B}$ only depend on $M$, $\hat M$, and $j^{k_0}_{p_0} f$.
\end{thm}

\begin{proof} Assuming that $\hat M$ is strictly pseudoconvex, we choose a
local basis of CR vector fields tangent to $\hat M$ such that
$\hat h_{\bar A B} (0) = i \da_A^B$. Since $f$ is transversal, it
follows that $M$ is strictly pseudoconvex at $0$, and we likewise
choose a local basis of CR vector tangent to $M$ such that
$h_{\bar a b} (0) = i \da_a^b$.
 As in
Theorem~\ref{T:strpcx}, we assume that $f$ is constantly
$(k_0,s)$-degenerate at $p_0$. We let
\[ V_k = \spanc \{ g_{\stringA B} (0) \colon \ell (\stringA) \leq
k \}
\]
and claim that we can assume that $V_{k_0} \subset \Cnh$ is just
the subspace spanned by the first $s$ unit vectors. To see this,
recall that
\[g_{\bar a_1\dots \bar a_k B} = \br{D_{\bar a_k} D_{\bar a_{k-1}}
\dots D_{\bar a_1} \tau , \hat L_B}.
\]
Hence, under a change of basis $\hat L_B^\p = u_B^C \hat L_C$,
the $g$ transform as $g_{\stringA B}^\p = u_B^D g_{\stringA D}$.
Also, since
\[ \hat h_{\bar A B} = \br{ d \hat \gth, \hat L_{\bar A} \wedge \hat L_B},\]
the $\hat h_{\bar A B}$ transform as $\hat h_{\bar A B}^\p =
\overline{ u_A^C } u_B^D \hat h_{\bar C D}$, and any $u_A^B$
which is unitary at $0$ will respect the normalization condition
$\hat h_{\bar A B} (0) = i \da_A^B$. So, our claim is proved.

Since $f$ is constantly $(k_0,s)$-degenerate, we can choose
$\stringA_{1,a},\dots,\stringA_{s,a}$ for $a\in\{1,\dots,n\}$ and
$\stringA_1,\dots,\stringA_s$, whose length does not exceed $k_0$,
with such that for any $a\in \{1,\dots,n\}$ the vectors
\[v_{a,j} = \left(g_{\stringA_{j,a}
B} -
  \sum_{\ell(\stringA^\p) < \ell(\stringA_{j,a})}
  g_{\stringA^\p B} \, v_{ \stringA_j}^{\stringA^\p}
    \left( \overline{L^J \gm_c^D} \right)\right),
    \quad j=1,\dots,s,
    \]
as well as
\[v_{j} = \left(g_{\stringA_j
B} -
  \sum_{\ell(\stringA^\p) < \ell(\stringA_j)}
  g_{\stringA^\p B} \, p_{a \stringA_j}^{\stringA^\p}
    \left( \overline{L^J \gm_c^D} \right)\right),
    \quad j=1,\dots,s,
    \]
 are linearly independent at $0$, and hence each of these $n+1$
 families of vectors span
$V_{k_0}$. For any  $t\in\{1,\dots,n\}$, $C\in\{1,\dots,\hat n\}$,
and $I=(I_1,\dots,I_k)$, the determinants of
\begin{equation}\label{E:detvanishesstrpscx}
  \begin{pmatrix}
    v_{a,1}^1 & \dots & v_{a,1}^s & v_{a,1}^C\\
    \vdots & & \vdots &\vdots  \\
    v_{a,s}^1 & \dots & v_{a, s}^s & v_{a,s}^C \\
    L^{\bar I} g_{\bar a 1} & \dots & L^{\bar I} g_{\bar a s} &L^{\bar I}
g_{\bar a C}
  \end{pmatrix}, \quad a=1,\dots,n,
\end{equation}
and
\begin{equation}\label{E:det2vanishesstrpscx}
  \begin{pmatrix}
    v_{1}^1 & \dots & v_{1}^s & v_{1}^C\\
    \vdots & & \vdots &\vdots  \\
    v_{s}^1 & \dots & v_{ s}^s & v_{s}^C \\
    L^{\bar I} g_{\bar t 1} & \dots & L^{\bar I} g_{\bar t s} &L^{\bar I}
g_{\bar t C}
  \end{pmatrix}
\end{equation}
vanish on $M$. Since
\[ L^{\bar I} g_{\bar t B} = \overline{ L^I\left( \gm_t^A \hat h_{A \bar
B}\right)} = \overline{ L^I \gm_t^A } \hat h_{\bar A B} +
q^I_{t,B} \left( \overline{L^J \gm_c^D}\right),
\]
where $|J|<k$, expanding \eqref{E:detvanishesstrpscx} and
\eqref{E:det2vanishesstrpscx} along the last row we obtain
\begin{equation}\label{E:detstrpscxexpanded}
  \Da^{\bar b}_{\bar a, C_0} \left( L^I \gm_a^A \right) \hat h_{A\bar b}
  + \bar \Da_a \left(  L^I \gm_a^A \right) \hat h_{A \bar C_0}=
  r^I_{a,C_0} \left({L^J \gm_a^B}\right), \quad
  a=1,\dots,n,
\end{equation}
\begin{equation}\label{E:det2strpscxexpanded}
  \Da^{\bar b}_{\bar C_0} \left(L^I \gm_t^A\right) \hat h_{A\bar b}
  + \bar \Da  \left(L^I \gm_t^A\right) \hat h_{A \bar C_0}=
  r^I_{C_0} \left({L^J \gm_a^B}\right), \quad
\end{equation}
where again $|J|<k$, $r^I_{a,C_0}$ and $r^I_{C_0}$ are polynomial
in their arguments with smooth coefficients, and $\Da^{\bar
b}_{\bar a, C_0}$, $\bar \Da_a$, $\Da^{\bar b}_{\bar C_0}$, $\bar
\Da_a$ denotes the various cofactors produced by the row
expansion. These are rational in $L^J \gm_a^B$, where $|J|\leq
k_0-1$, and polynomial in $L^K \xi$, where $|K|\leq k_0$.

We need the following computational fact, which the reader can
verify easily.

\begin{lem}\label{L:interchanging1}
Given $P = p_1 \dots p_k$ with $1\leq p_j\leq n$, the following
holds for any $a$.
\begin{equation}\label{E:interchanging1}
  L_{\bar a} L^P = L^P L_{\bar a} - \sum_{j=1}^k h_{\bar a p_j}
  L^{P_j} T + \sum_{|K| \leq k-2} C_K L^K T,
\end{equation}
where $P_j$ is obtained from $P$ by removing the $j$th entry.
\end{lem}

Now given any $I = i_1 \dots i_k$ with $i_1\leq i_2 \leq \dots
\leq i_k$, we first apply the operator $L^I$ to the equations
\eqref{E:complemma1} and \eqref{E:complemma2} to obtain
\begin{multline} \label{E:firstforetas} g_{\stringA B}\left( L^I \eta^B
\right) -
  \sum_{\ell(\stringA^\p) < \ell(\stringA) } g_{\stringA^\p B} \left(L^I
\eta^B\right) \, v_{\stringA}^{\stringA^\p}
    \left( \overline{L^J \gm_c^D}\right)\\ = \td w_{\stringA}
    \left(L^M \gm_c^D, \overline{L^J T^m \gm_c^D},\overline{L^J T^m \eta^D}
\right)
\end{multline}
and
\begin{multline}\label{E:thenforgammas}
     g_{\stringA B} \left( L^I \gm_a^B \right) -
  \sum_{\ell(\stringA^\p) < \ell(\stringA)}
  g_{\stringA^\p B} \left( L^I \gm_a^B \right) \, p_{a \stringA}^{\stringA^\p}
    \left( \overline{L^J \gm_c^D}\right) - \\
  \sum_{\ell(\stringA^\p) < \ell(\stringA) } g_{\stringA^\p B} \left( L^I
\eta^B \right) \, q_{a \stringA}^{\stringA^\p}
    \left( \overline{L^J \gm_c^D}\right) = \td r_{a \stringA} \left(L^M
\gm_c^D, \overline{L^J T^m \gm_c^D},
    \overline{L^J T^m \eta^D}\right),
\end{multline}
where $|M|<k$,  $|K|\leq k_0 -1$, $|J|\leq k_0$, $m\leq k$, and
$|J| + m \leq k_0$. We have also simplified matters a bit by
applying \eqref{E:levitransform} in order to get rid of the terms
involving $\xi$.

Let us recall \eqref{E:det2strpscxexpanded}:
\[ \Da^{\bar b}_{\bar C_0} \left(L^{I} \gm_t^A\right) \hat h_{A\bar b}
  + \bar \Da  \left(L^{I} \gm_t^A\right) \hat h_{A \bar C_0}=
  r^I_{C_0} \left({L^M \gm_a^B}\right),\]
where  $|M| \leq k-1$. Applying $L_{\bar t}$ to this equation, and
applying Lemma \ref{L:interchanging1} and
\eqref{E:Tderivativeofgamma}, we obtain
\begin{multline}\label{E:missingequationforeta}
\left(L^{I} \gm_t^A\right) L_{\bar t} \left( \Da^{\bar b}_{\bar
C_0}  \hat h_{A\bar b}\right)
  + \left(L^{I} \gm_t^A\right) L_{\bar t} \left( \bar \Da   \hat h_{A \bar
C_0}
  \right)
  - \Da^{\bar b}_{\bar C_0} \left(L^{I} \eta^A\right) h_{\bar t t} \hat
h_{A\bar
  b}\\
  - \bar \Da  \left(L^{I} \eta^A\right) h_{\bar t t} \hat h_{A \bar
  C_0}
  - \sum_{j=1}^k h_{\bar t i_j} \left( L^{I_j} L_t \eta^A \right) \Da^{\bar
b}_{\bar C_0}\hat h_{A\bar b}
  - \sum_{j=1}^k h_{\bar t i_j} \left( L^{I_j} L_t \eta^A \right) \bar \Da
\hat h_{A\bar C_0}
  \\
  =
  \td r^I_{C_0} \left(L^M \eta^B, {L^M \gm_a^B}, L^J T \eta^B, {L^J T
\gm_a^B}\right).
\end{multline}
In this equation, $|M|<k$, and $|J|<k-1$.

We consider the system of linear equations in the $L^I \gm_a^B$
and $L^I \eta^B$ where $I$ ranges over all $I = i_1 \dots i_k$
with $i_1 \leq i_2 \leq \dots \leq i_k$, $a\in \{ 1,\dots , n\}$,
$B\in \{1\dots,\hat n \}$ which is obtained in the following way.
In \eqref{E:firstforetas} we set $\stringA = \stringA_1 ,\dots ,
\stringA_s$. For each $I$, this gives $s$ linear equations in $L^I
\eta^B$. We next consider the $\hat n - s$ equations obtained from
\eqref{E:missingequationforeta} by replacing $t$ by $i_k$ and
letting $C_0 = s+1,s+2,\dots,\hat n$. For the equations in $L^I
\gm_a^B$ ($a\in\{1,\dots ,n\}$  fixed) we replace $\stringA$ in
\eqref{E:thenforgammas} by $\stringA_{a,1}\dots \stringA_{a,s}$,
which gives us (for each $I$ and each $a$) $s$ linear equations in
$L^I \gm_a^B$. For each $I$ and each $a$ we also consider the
equations obtained from \eqref{E:det2strpscxexpanded} by replacing
$t$ by $a$, and letting $C_0 = s+1,\dots ,\hat n$.

With the simplifications above, it is easy to see that this linear
system can be solved for $L^I \gm_a^B$ and $L^I \eta^B$. We can
now lift the restriction that $I$ is increasing by noting that
rearranging the order of the $i_j$s only produces error terms of
order less than $k$. This finishes the proof.
\end{proof}

\section{Constructing the complete system; proof of Theorem
\ref{T:complsys}}\label{S:completesystem} We are now going to
construct the complete system for the jets of $f$. Since Theorems
\ref{T:reflectionlevinondegenerate} and \ref{T:reflectionstrpscx}
are a bit different than the reflection identities in \cite{E5},
we will go into some detail. For the rest of this section we
assume that we are either in the situation of Theorem
\ref{T:reflectionlevinondegenerate} or in the situation of Theorem
\ref{T:reflectionstrpscx}. First, let us recall the following fact
from \cite{E5} (the second part of Proposition 3.18):

For any multi-index $J$, integer $k\geq 1$, and index
$c\in\{1,\ldots, n\}$ there exist smooth functions  $b^{e_1\ldots
e_m}_s$  such that
\begin{equation}\label{E:maincommutatoridentity}
\sum_{m=1}^{|J|+k}\sum_{s=0}^{k} b_s^{e_1\ldots e_m}
\underbrace{[\ldots [L_{e_1}\ldots L_{e_m},L_{\bar c}],L_{\bar
c}]\ldots,L_{\bar c}]}_{\text{\rm length $s$ }}=(h_{\bar c
1})^{p}L^{J}T^k.
\end{equation}
Here, $p=k+|J|-|J|_1+1$, $|J|_1$ denotes the number of occurences
of the index $1$ in the multi-index $J$, and the length of the
commutator $[\ldots [X,Y_{1}],Y_{2}]\ldots,Y_{s}]$ is $s$. Since
in any case $M$ is Levi-nondegenerate at $0$, we can use this
identity in order to find a formula for $L^J T^k \gm_a^B$ and $L^J
T^k \eta^B$. We first apply $L_{\bar c}$ (at most $k$ times) to
$\gm_a^B$ and $\eta_a^B$, and using equations
\eqref{E:barredderivativeofgamma} and
\eqref{E:barredderivativeofeta} we see that the result is
polynomial in $\gm_c^D$, $\eta_c^D$, and $\overline{L_{\bar c}^j
\gm_c^D}$, where $j\leq k-1$. Next, we apply $L_{e_1}\dots
L_{e_m}$ to the result, and we apply our reflection identities
(either \eqref{E:reflectionlevinondegenerategamma} and
\eqref{E:reflectionlevinondegenerateeta} or
\eqref{E:reflectionpseudoconvexgamma} and
\eqref{E:reflectionpseudoconvexeta}) to see that the result is
rational in
\[L^M T^l \gm_c^D, L^M T^l \eta^D,
 \overline{ L^K T^m \gm_c^D}, \overline{ L^K T^m \eta^D},
\]
where $|M|+l \leq \max(k_0,|J| + k -1)$,  $|K| + m  \leq k_0 + 1$,
and $m\leq |J|$. Another application of $L_{\bar c}$ (at most $k$
times) then gives us that
\begin{align}\label{E:LTgammaisequal1}
L^J T^k \gm_a^B &= u(L^M T^l \gm_c^D, L^M T^l \eta^D,
 \overline{ L^K T^m \gm_c^D}, \overline{ L^K T^m \eta^D}), \\
 \label{E:LTetaisequal1}
 L^J T^k \eta^B &= v(L^M T^l \gm_c^D, L^M T^l \eta^D,
 \overline{ L^K T^m \gm_c^D}, \overline{ L^K T^m \eta^D}),
\end{align}
where $|M|+l \leq \max(k_0,|J| + k -1)$,  $|K| + m  \leq k + k_0 +
1$, and $m\leq \min(|J| + k, k_0)$. We are abusing notation by
denoting by $u$, $v$ and $w$ functions which are rational in their
arguments with coefficients smooth on $M\times\hat M$ and which
only depend on $M$ and $\hat M$ (and are allowed to change). Now
in these equations for $|J| + k = 2 k_0 + 2$ we can substitute for
the conjugated terms, yielding the same equations, but with the
bounds $|M|+l \leq 2k_0 + 1 $, $|K| + m \leq 2k_0 + 1$. Another
use of equations \eqref{E:barredderivativeofgamma} and
\eqref{E:barredderivativeofeta} lets us conclude that for all $J$,
$k$, $R$ with $|J| + k + |R| = 2k_0 + 2$,
\begin{align}\label{E:LTLbargammaisequal}
L^J T^k L^{\bar R}\gm_a^B &= u(L^M T^l \gm_c^D, L^M T^l \eta^D,
 \overline{ L^K T^m \gm_c^D}, \overline{ L^K T^m \eta^D}), \\
 \label{E:LTLbaretaisequal1}
 L^J T^k L^{\bar R}\eta^B &= v(L^M T^l \gm_c^D, L^M T^l \eta^D,
 \overline{ L^K T^m \gm_c^D}, \overline{ L^K T^m \eta^D}),
\end{align}
with $|M|+l \leq 2k_0 + 1 $, $|K| + m  \leq 2k_0 + 1$. These two
equations together with \eqref{E:levitransform} let us also see
that for all $J$, $k$, $R$ with $|J| + k + |R| = 2k_0 + 2$,
\begin{equation}\label{E:LTLxiisequal}
  L^J T^k L^{\bar R}\xi = w(L^M T^l \gm_c^D, L^M T^l \eta^D,
 \overline{ L^K T^m \gm_c^D}, \overline{ L^K T^m \eta^D}),
\end{equation}
again with $|M|+l \leq 2k_0 + 1 $, $|K| + m  \leq 2k_0 + 1$.
Equations \eqref{E:LTLbargammaisequal},
\eqref{E:LTLbaretaisequal1}, and \eqref{E:LTLxiisequal}, written
in some local coordinate systems, constitute the complete system
of differential equations described in Theorem \ref{T:complsys}.
Also note that the same complete system works for mappings whose
jets are close to the $2k_0 + 2$ jet of $f$. This finishes the
proof of Theorem \ref{T:complsys}.

\bibliographystyle{plain}
\bibliography{050102}
\end{document}